\documentclass[journal]{IEEEtran}

\usepackage{color}
\usepackage{soul}
\usepackage{subcaption}

\usepackage{amssymb, amsmath} 
\usepackage{graphicx}
\usepackage{cite}
\usepackage{epstopdf}
\usepackage[margin=2cm]{geometry}
\usepackage{latexsym}
\usepackage{epsfig}
\usepackage{mathrsfs}
\usepackage{bbm}
\usepackage{color}

\newtheorem{thm}{\bf{Theorem}}[section]

\newtheorem{proposition}{\bf{Proposition}}[section]

\newtheorem{remark}{\bf{Remark}}[section]

\newcommand{\nn}{\nonumber}
\newcommand{\vu}{\vec{u}}
\newcommand{\vm}{\vec{m}}
\newcommand{\vx}{\vec{x}}
\newcommand{\vw}{\vec{w}}
\newcommand{\vy}{\vec{y}}

\newcommand{\vM}{\vec{M}}
\newcommand{\vX}{\vec{X}}
\newcommand{\vW}{\vec{W}}
\newcommand{\vY}{\vec{Y}}

\newcommand{\vb}{\vec{b}}

\newcommand{\vC}{\vec{C}}
\newcommand{\vL}{\vec{L}}
\newcommand{\ve}{\vec{e}}

\DeclareMathOperator{\tr}{tr}

\begin{document}
\sloppy
\title{Quadratic Multi-Dimensional Signaling Games and Affine Equilibria
\thanks{S. Sar{\i}ta\c{s} and S. Gezici are with the Department of Electrical and Electronics Engineering, Bilkent University, 06800, Ankara, Turkey. Emails: \{serkan,gezici\}@ee.bilkent.edu.tr. S. Y\"uksel is with the Department of Mathematics and Statistics, Queen's University, Kingston, Ontario, Canada, K7L 3N6.  Email: yuksel@mast.queensu.ca.}
\thanks{This research was supported in part by the Natural Sciences and Engineering Research Council (NSERC) of Canada, The Scientific and Technological Research Council of Turkey (T\"{U}B\.{I}TAK) and the Distinguished Young Scientist Award of Turkish Academy of Sciences (T\"{U}BA-GEB\.{I}P 2013).}
\thanks{Part of this work was presented at the 2015 American Control Conference (ACC), Chicago, IL, 2015.}
}
\author{Serkan Sar{\i}ta\c{s}, Serdar Y\"uksel, and Sinan Gezici}
\maketitle
\begin{abstract}
This paper studies the decentralized quadratic cheap talk and
signaling game problems when an encoder and a decoder, viewed as two
decision makers, have misaligned objective functions. The main
contributions of this study are the extension of Crawford and Sobel's
cheap talk formulation to multi-dimensional sources and to noisy
channel setups. We consider both (simultaneous) Nash equilibria and (sequential)
Stackelberg equilibria. We show that for arbitrary scalar sources, in the presence of misalignment, the quantized nature of all equilibrium
policies holds for Nash equilibria in the sense that all Nash
equilibria are equivalent to those achieved by quantized encoder
policies. On the other hand, all Stackelberg equilibria policies are
fully informative. For multi-dimensional setups, unlike the scalar
case, Nash equilibrium policies may be of non-quantized nature, and
even linear. In the noisy setup, a Gaussian source is to be
transmitted over an additive Gaussian channel. The goals of the
encoder and the decoder are misaligned by a bias term and encoder's
cost also includes a penalty term on signal power. Conditions for the
existence of affine Nash equilibria as well as general informative
equilibria are presented. For the noisy setup, the only Stackelberg
equilibrium is the linear equilibrium when the variables are scalar.
Our findings provide further conditions on when affine policies may
be optimal in decentralized multi-criteria control problems and lead
to conditions for the presence of active information transmission in
strategic environments. 
\end{abstract}


\section{Introduction}

Team theory is concerned with the interaction dynamics among decentralized decision makers with identical objective functions. On the other hand, game theory deals with setups with misaligned objective functions, where each player chooses a strategy to maximize its own utility which is determined by the joint strategies chosen by all players. Information transmission in team problems is well-understood with extensive publications present in the literature; for a detailed account we refer the reader to \cite{YukselBasarBook}. Despite the difficulty to obtain solutions under general information structures, it is evident in team problems that more information provided to any of the decision makers does not hurt the system performance and there is a well-defined partial order of information structures as studied by Blackwell \cite{Blackwell} and others. However, for general non-zero sum game problems, informational aspects are very challenging to address; more information can hurt some or even all of the players in a system, see e.g. \cite{hirshleifer1971private}. Further intricacies on informational aspects in competitive setups have been discussed in \cite{gossner2001value}, \cite{kamien1990value} and\cite{tBasarStochasticDiffGames}.

Signaling games and cheap talk are concerned with a class of Bayesian games where an informed decision maker transmits information to another decision maker. Unlike a team setup, however, the goals of the agents are misaligned. Such a study has been initiated by Crawford and Sobel \cite{SignalingGames}, who obtained the striking result that under some technical conditions on the utility functions of the decision makers, the cheap talk problem only admits equilibrium policies that are essentially quantization policies. This is in significant contrast with the case where the utility functions are aligned. 

The cheap talk and signaling game problems find applications in networked control systems when a communication channel/network is present among competitive and non-cooperative decision makers \cite{basols99}. For example, in a smart grid application, there may be strategic sensors in the system \cite{misBehavingAgents} that wish to alter the equilibrium decisions at a controller receiving data from the sensors to lead to a more desirable equilibrium, for example by enforcing an outcome to enhance its prolonged use in the system. One may also consider a utility company which wishes to inform users regarding pricing information; if the utility company and the users engage in selfish behaviour, it may be beneficial for the utility company to hide certain information and the users to be strategic about how they interpret the given information. One further area of application is recommender systems (as in rating agencies) \cite{miklos2013value}. All of these applications lead to a drastically new framework where the value of information and its utilization are very fragile to the system under consideration and our study here is an initiator for such a general setup. 

Even though in this paper we only consider quadratic criteria under a bias term leading to a misalignment, the contrast with the case where there is no bias (that has been heavily studied in the information theory literature) raises a number of sharp conclusions for system designers working on networked systems under competitive environments.

Identifying when optimal policies are linear or affine for decentralized systems involving Gaussian variables under quadratic criteria is a recurring problem in control theory, starting perhaps from the seminal work of Witsenhausen \cite{wit68}, where suboptimality of linear policies for such problems under {\it non-classical information structures} is presented. The reader is referred to {\sl Chapters 3} and {\sl 11} of \cite{YukselBasarBook} for a detailed discussion on when affine policies are and are not optimal. These include the problem of communicating a Gaussian source over a Gaussian channel, variations of Witsenhausen's counterexample  \cite{banbas87a}; and game theoretic variations of such problems. For example if the noise variable is viewed as the maximizer and the encoders/decoders (or the controllers) act as the minimizer, then affine policies may be optimal for a class of settings, see  \cite{basmin71,basmin72,RotkowitzUniform,GattamiRobust,basCDC2008}. \cite{basCDC2008} also provides a review on Linear Quadratic Gaussian (LQG) problems under nonclassical information including Witsenhausen's counterexample. Our study provides further conditions on when affine policies may constitute equilibria for such decentralized quadratic Gaussian optimization problems.

There have been a number of related contributions in the economics literature in addition to the seminal work by Crawford and Sobel, which we briefly review in the following: Reference \cite{complexityConflict} shows that even if the sender and the receiver have identical preferences, perfect communication may not be possible in an equilibrium because information transmission may be costly. Reference \cite{Krishna2001} studies the setup in \cite{SignalingGames} with two senders and shows that if senders transmit the messages sequentially once, then the equilibrium is always quantized and if senders transmit the messages simultaneously and their biases are either both positive or both negative, then a fully revealed equilibrium is possible. Reference \cite{shintaromiura2012} studies a scalar setup and proves that if multiple senders transmit the messages sequentially and their biases have opposite signs, then a fully revealed equilibrium is possible; this study also considers two-dimensional real valued sources, and shows that a fully revealed equilibrium occurs if and only if the multiple senders have perfectly opposing biases. Moreover, multi-dimensional cheap talk with multiple senders is analyzed in \cite{BattagliniMultiCheapTalk} and \cite{multiRestricted} with unbounded and bounded state spaces, respectively. The study in \cite{blumeNoisyTalk} considers a special noisy channel setup between the sender and receiver, and shows that there may be infinitely many actions (countable or uncountable) induced in an equilibrium even though all equilibria are interval partitions in the noiseless case \cite{SignalingGames}. Conditions for Nash equilibria are investigated in \cite{ncnc} for a scenario in which there exists a discrete noisy channel between an informed sender and an uninformed receiver, and the source is finitely valued. Furthermore, there are some contributions which modify the information structure given in Crawford and Sobel's setup: In \cite{chenying2009}, the sender knows that the receiver has partial information about his/her private information; whereas the sender does not know this in \cite{twoSidedCheap, expertAdviceTwoSidedCheap}. Reference \cite{dynamicStrategicInfoTrans} studies Crawford and Sobel's setup in a finite horizon environment where, in each period, a privately informed sender transmits a message and a receiver takes an action. For a detailed literature review on communication between informed experts and uninformed decision makers, we refer the reader to \cite{givingReceivingAdvice}. We note also that in the area of information theory, there exists a vast literature on security aspects of information transmission, see e.g., \cite{yamamoto1988rate,tandon2013}. Game theoretic analysis is also useful in various contexts involving security problems. For example, the security of the smart-grid infrastructure can be analyzed by considering the adversarial nature of the interaction between an attacker and a defender \cite{cyberSecuritySmartGrid, networkSecurity}, and a game theoretic setup would be appropriate to analyze such interactions. For an overview of security and privacy problems in computer networks that are analyzed within a game-theoretic framework, \cite{gameTheoryNetworkSecurity} can be referred.

In the control community, recently, there have been few studies: \cite{CedricWork} considered a Gaussian cheap talk game with quadratic cost functions where the analysis considers Stackelberg equilibria, for a class of single- and multi-terminal setups and where linear equilibria have been studied. For the setup of Crawford and Sobel but when the source admits an exponentially distributed real random variable, \cite{FurrerReport} establishes the discrete-nature of equilibria, and obtains the equilibrium bins with {\it finite upper bounds on the number of bins} under any equilibrium in addition to some structural results on informative equilibria for general sources. 

\subsection{Contributions} 
The main contributions of this study are as follows. We prove that for any scalar source, all Nash equilibrium policies at the encoder are equivalent to some quantized policy, but all Stackelberg equilibrium policies are fully informative. That is, there is some information hiding for the Nash setup, as opposed to the Stackelberg setup. We show that for multi-dimensional setups, however, unlike the scalar case, Nash equilibrium policies may be non-quantized and can in fact be linear. In the noisy setup, a Gaussian source is to be transmitted over an additive Gaussian channel. The goals of the encoder and the decoder are misaligned by a bias term and encoder's cost also includes a penalty term of the transmitted signal. Conditions for the existence of affine equilibrium policies as well as general informative Nash equilibria are presented for both the scalar and multi-dimensional setups. We compare the results with socially optimal costs and information theoretic lower bounds, and discuss the effects of the bias term on equilibria. Furthermore, we prove that the only equilibrium in the Stackelberg noisy setup is the linear equilibrium for the scalar case.

\section{Problem Definition}

\begin{figure} [ht]
	\begin{center}
		\includegraphics[scale=0.5]{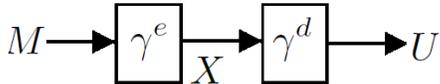}
	\end{center}
	\caption{General system model for noiseless case.}
	\label{figure:noiselessScheme}
\end{figure}

Let there be two decision makers (DMs): An encoder (DM 1) and a decoder (DM 2) as shown in Fig.~\ref{figure:noiselessScheme}. DM 1 wishes to encode the $\mathbb{M}$-valued random variable $M$ to DM 2. Let $X$ denote the $\mathbb{X}$-valued random variable which is transmitted to DM 2. DM 2, upon receiving $X$, generates its optimal decision $U$ which we also take to be $\mathbb{M}$-valued. We allow for randomized decisions, therefore, we let the policy space of DM 1 be the set of all stochastic kernels from $\mathbb{M}$ to $\mathbb{X}$.\footnote{Recall that $P$ is a stochastic kernel from $\mathbb{M}$ to $\mathbb{X}$ if $P(\cdot|m)$ is a probability measure on ${\cal B}(\mathbb{X})$ for every $m\in \mathbb{M}$ and for every Borel $A\in{\cal B}(\mathbb{X})$, $P(A|\cdot)$ is a Borel measurable function of $m$.} Let $\Gamma^e$ denote the set of all such policies. We let the policy space of DM 2 be the set of all stochastic kernels from $\mathbb{X}$ to $\mathbb{M}$. Let $\Gamma^d$ denote the set of all such stochastic kernels.

Given $\gamma^e \in \Gamma^e$ and $\gamma^d \in \Gamma^d$, the goal in the classical communications theory is to minimize the expectation
\[J(\gamma^e,\gamma^d) = \int c(m, u) \gamma^e(\mathrm{d}x|m)\gamma^d(\mathrm{d}u|x) P(\mathrm{d}m),\]
where $c$ is some cost function. One very common case is the setup with $c(m,u)=|m-u|^2$. 

Recall that a collection of decision makers who have an agreement on the probabilistic description of a system and a cost function to be minimized, but who may have different on-line information is said to be a {\it team} (see, e.g. \cite{YukselBasarBook}). Hence, the classical communications setup may be viewed as a team of an encoder and a decoder.

In many applications (in networked systems, recommendation systems, and applications in economics) the objectives of the encoder and the decoder may not be aligned. For example, DM 1 may aim to minimize
\begin{eqnarray}
J^e(\gamma^e,\gamma^d) = \int c^e(m, u) \gamma^e(\mathrm{d}x|m)\gamma^d(\mathrm{d}u|x) P(\mathrm{d}m)\,,\nn
\end{eqnarray}
whereas DM 2 may aim to minimize
\begin{eqnarray}
J^d(\gamma^e,\gamma^d) = \int c^d(m, u) \gamma^e(\mathrm{d}x|m)\gamma^d(\mathrm{d}u|x) P(\mathrm{d}m)\,.\nn
\end{eqnarray}
In this study, the problems are investigated where the encoder and the decoder are deterministic rather than randomized; i.e., $\gamma^e(\mathrm{d}x|m)=\mathbbm{1}_{\{f^e(m)\in \mathrm{d}x\}}$ and $\gamma^d(\mathrm{d}u|x)=\mathbbm{1}_{\{f^d(x)\in \mathrm{d}u\}}$ where $\mathbbm{1}_{\{D\}}$ denotes the indicator function of an event $D$, and $f^e(m)$ and $f^d(x)$ are some deterministic functions of the encoder and decoder, respectively. Such a problem is known in the economics literature as {\it cheap talk} (the transmitted signal does not affect the cost, that is why the game is named as {\it cheap talk}). A more general formulation would be the case when the transmitted signal is also an explicit part of the cost function $c^e$ or $c^d$; in that case, the setup is called a {\it signaling game}. We will consider a noisy communication setup, where the problem may be viewed as a signaling game, rather than cheap talk, later in this study.

Since the goals are not aligned, such a problem is studied under the tools and concepts provided by {\it game theory}. A pair of policies $\gamma^{*,e}, \gamma^{*,d}$ is said to be a {\bf Nash equilibrium} if
\begin{align*}
	\begin{split}
	 J^e(\gamma^{*,e}, \gamma^{*,d}) &\leq J^e(\gamma^{e}, \gamma^{*,d}) \quad \forall \gamma^e \in \Gamma^e \,, \\
 J^d(\gamma^{*,e}, \gamma^{*,d}) &\leq J^d(\gamma^{*,e}, \gamma^{d}) \quad \forall \gamma^d \in \Gamma^d \,.
	\end{split}
	\end{align*}
We note that when $c^e=c^d$ the setup is a traditional communication theoretic setup. If $c^e = -c^d$, that is, if the setup is a zero-sum game, then an equilibrium is achieved when $\gamma^e$ is non-informative (e.g., a kernel with actions statistically independent of the source) and $\gamma^d$ uses only the prior information (since the received information is non-informative). We call such an equilibrium a {\it non-informative (babbling) equilibrium}. The following is a useful observation, which follows from \cite{SignalingGames}:

\begin{proposition} \label{prop:nonInEq}
	A non-informative (babbling) equilibrium always exists for the cheap talk game.
\end{proposition}

In the discussion so far, a simultaneous game-play is assumed and thus equilibrium refers to a Nash equilibrium. Besides the simultaneous game-play, one can also consider a sequential game-play; i.e. first the encoder sends the message, then the decoder receives it and takes an action sequentially while first the encoder's policy is announced. {\bf Stackelberg} equilibria arise in this case. In the Stackelberg game, the encoder announces his coding strategy and since the decoder takes an action after receiving the message, the encoder knows the optimal action which will be taken by the decoder and chooses the message to be transmitted accordingly. A pair of policies $\gamma^{*,e}, \gamma^{*,d}$ is said to be a {\bf Stackelberg equilibrium} if
\[ J^e(\gamma^{*,e}, \gamma^{*,d}(\gamma^{*,e})) \leq J^e(\gamma^e, \gamma^{*,d}(\gamma^e)) \quad \forall \gamma^e \in \Gamma^e\,, \]
where $\gamma^{*,d}(\gamma^e)$ satisfies
\[ J^d(\gamma^{e}, \gamma^{*,d}(\gamma^{e})) \leq J^d(\gamma^{e}, \gamma^d(\gamma^{e})) \quad \forall \gamma^d \in \Gamma^d.\]
Throughout the paper, all equilibrium terms refer to the Nash equilibrium unless otherwise stated; it will be separately indicated for the Stackelberg game setup and equilibrium.

Crawford and Sobel \cite{SignalingGames} have made foundational contributions to the study of cheap talk with misaligned objectives where the cost functions $c^e$ and $c^d$ satisfy certain monotonicity and differentiability properties but there is a bias term in the cost functions. Their result is that the number of bins in an equilibrium is upper bounded by a function which is negatively correlated to the bias.

We will first consider the scalar setting by taking the cost functions as $c^e\left(m,u\right) = \left(m-u-b\right)^2$ and $c^d\left(m,u\right) = \left(m-u\right)^2$ where $b$ denotes the bias term. The motivation for such functions stems from the fields of information theory, communication theory and LQG control; for these fields quadratic criteria are extremely important. Recall that for the case with $b=0$, the cost functions simply reduce to those for a minimum mean-square estimation (MMSE) problem. 

\section{Quadratic Cheap Talk}

\subsection{Nash Equilibria in the Scalar Case}
As before, let the cost functions be defined as \mbox{$c^d\left(m,u\right) = \left(m-u\right)^2$} and \mbox{$c^e\left(m,u\right) = \left(m-u-b\right)^2$} where $b$ is the bias term. Some existence and deterministic properties of the equilibrium policies of the encoder and the decoder are stated in \cite{FurrerReport} and \cite[Chp.4]{YukselBasarBook}.
\begin{thm}\cite{FurrerReport}
(i) For any $\gamma^e$, there exists an optimal $\gamma^d$, which is deterministic. (ii)
For any $\gamma^d$, any randomized encoding policy can be replaced with a deterministic $\gamma^e$ without any loss to DM 1. (iii) Suppose $\gamma^e$ is an $M$-cell quantizer, then there exists an optimal $\gamma^d$, which is the conditional expectation of the respective bin.
\end{thm}

The following builds on \cite[Lem.1]{SignalingGames}, which considers sources on $[0,1]$ that admit densities. We note that the analysis here applies to arbitrary scalar valued random variables. The proofs essentially follow from \cite{SignalingGames}.
\begin{thm} \label{noiselessDiscrete}
Let $m$ be a real-valued random variable with an arbitrary probability measure. Let the strategy set of the encoder (DM 1) consists of the set of all measurable (deterministic) functions from $\mathbb{M}$ to $\mathbb{X}$. Then, an equilibrium encoder policy has to be quantized almost surely, that is, it is equivalent to a quantized policy for the encoder in the sense that the performance of any equilibrium encoder policy is equivalent to the performance of a quantized encoder policy. Furthermore, the quantization bins are convex.
\end{thm}
\begin{IEEEproof}
Let there be an equilibrium in the game (with possibly uncountably infinitely many bins, countably many bins or finitely many bins). Let two bins be $\mathcal{B}^{\alpha}$ and $\mathcal{B}^{\beta}$. Also let $m^{\alpha}$ indicate any point in $\mathcal{B}^{\alpha}$; i.e., $m^{\alpha} \in \mathcal{B}^{\alpha}$. Similarly, let $m^{\beta}$ represent any point in $\mathcal{B}^{\beta}$; i.e., $m^{\beta} \in \mathcal{B}^{\beta}$. The decoder chooses action $u^{\alpha}=\mathbb{E}[m|m \in \mathcal{B}^{\alpha}]$ when the encoder sends $m^{\alpha} \in \mathcal{B}^{\alpha}$ and action $u^{\beta}=\mathbb{E}[m|m \in \mathcal{B}^{\beta}]$ when the encoder sends $m^{\beta} \in \mathcal{B}^{\beta}$ in order to minimize its total cost. Without loss of generality, we can assume that  $u^{\alpha}<u^{\beta}$. Let $F(m,u)\triangleq(m-u-b)^2$. Because of the equilibrium definitions from the view of the encoder; $F(m^{\alpha},u^{\alpha})<F(m^{\alpha},u^{\beta})$ and $F(m^{\beta},u^{\beta})<F(m^{\beta},u^{\alpha})$. Hence, $\exists$ $\overline{m}$ that satisfies $F(\overline{m},u^{\alpha})=F(\overline{m},u^{\beta})$ which reduces to 
\begin{align}
\overline{m}=\frac{u^{\alpha}+u^{\beta}}{2}+b \iff (\overline{m}-u^{\alpha}) = (u^{\beta}-\overline{m}) +2b
\label{eq:mBarEquality}
\end{align}
Since $F(\overline{m}+\Delta,u^{\alpha})>F(\overline{m}+\Delta,u^{\beta})$ for any $\Delta>0$, $\mathcal{B}^{\beta}$ and $\{m|m<\overline{m}\}$ are disjoint sets. Similarly, $\mathcal{B}^{\alpha}$ and $\{m|m>\overline{m}\}$ are disjoint sets, too. Thus, from the definitions of $u^{\alpha}$ and $u^{\beta}$, we have $u^{\alpha}<\overline{m}<u^{\beta}$ which implies $\overline{m}-u^{\alpha}>0$ and $u^{\beta}-\overline{m}>0$. Then, from \eqref{eq:mBarEquality},
\begin{align*}
u^{\beta} - u^{\alpha} &=  (u^{\beta}-\overline{m}) + (\overline{m}-u^{\alpha}) = 2(u^{\beta}-\overline{m})+2b\\
&>2b
\end{align*}
and
\begin{align*}
u^{\beta} - u^{\alpha} &=  (u^{\beta}-\overline{m}) + (\overline{m}-u^{\alpha}) = 2(\overline{m}-u^{\alpha})-2b\\
&>-2b
\end{align*}
are obtained. Hence, $u^{\beta}-u^{\alpha}>2|b|$, which implies that there must be at least $2|b|$ distance between the equilibrium points (decoder's actions, centroids of the bins). Further, from the encoder's point of view, given any two bins $\mathcal{B}^{\alpha}$ and $\mathcal{B}^{\beta}$, there exists a point $\overline{m}$ which lies between these two bins. This assures that each bin must be a single interval; i.e., convex cell except for a possible insignificant set of points with measure zero. Since there is an injective and monotonic relation between the convex cells of the encoder and decoder's actions, the equilibrium policy must be quantized almost surely.
\end{IEEEproof}

Recall again that for the case when the source admits density on $[0,1]$, Crawford and Sobel established the discrete nature of the equilibrium policies. For the case when the source is exponential, \cite{FurrerReport} established the discrete-nature, and obtained the equilibrium bins with {\it finite upper bounds on the number of bins} in any equilibrium. 

\subsection{Stackelberg Equilibria in the Scalar Case}
We will now observe that the Stackelberg setup is less interesting.
\begin{thm}
The Stackelberg equilibrium is unique and corresponds to a fully revealing (fully informative) encoder policy.
\label{thm:cheapStackelberg}
\end{thm}
\begin{IEEEproof}
Due to the Stackelberg assumption, the encoder knows that the decoder will use $\gamma^d(x) = u = \mathbb{E}[m|x]$ as an optimal decoder policy to minimize its cost. Then the goal of the encoder is to minimize the following:
\begin{align*}
\min_{x = \gamma^e(m)} & \mathbb{E} [(m-u-b)^2] = \min_{x = \gamma^e(m)} \mathbb{E} [(m-\mathbb{E}[m|x]-b)^2] \nn\\
&\qquad\qquad\qquad= \min_{x = \gamma^e(m)} \mathbb{E} [(m-\mathbb{E}[m|x])^2] + b^2 \nn\\
&\qquad\qquad\qquad= \min_{x = \gamma^e(m)} \mathbb{E} [(m-u)^2] + b^2 \,.
\end{align*}
Here, the second equality follows from the law of the iterated expectations. Since the goal of the decoder is to minimize $\min_{u = \gamma^d(x)} \mathbb{E} [(m-u)^2]$, the goals of the encoder and the decoder become essentially the same in the Stackelberg game setup, which effectively reduces the game setup to a {\it team} setup. In the team setup, the equilibrium is fully informative; i.e. the encoder reveals all of its information.\end{IEEEproof}

\subsection{Multi-Dimensional Cheap Talk: Nash Equilibria}

Our goal in this subsection is to show that it is possible to have linear equilibria in a multi-dimensional quadratic cheap talk, unlike the scalar setup. Let the source be uniform on $[0,1] \times [0,1]$ and the cost function of the encoder be defined by $c^e(\vm,\vu) = \|\vm-\vu-\vb\|^2$ and the cost function of the decoder be defined by $c^d(\vm,\vu) = \|\vm-\vu\|^2$ where the lengths of the vectors are defined in $L_2$ norm and $\vb$ is the bias vector. For such a scenario, we have the following result.

\begin{thm}\label{maxBin2D}
An equilibrium policy can be non-discrete and even linear.
\end{thm}
\begin{IEEEproof}
It suffices to provide an example. Consider $\vb = [0.3 \; 0]$. Then, as a (properly interpreted) limit case of the equilibrium in Fig.~\ref{figure:equilibriumInf2}, the following encoder and decoder policies form an equilibrium:
\begin{align*}
\gamma^e(m_1,m_2) &= (x_1,x_2) = (0, m_2)\,, \\
\gamma^d(x_1,x_2) &= (u_1,u_2) = (0.5, m_2)\,. 
\end{align*}
Here, the scalar setup is applied on the $x$-dimension with one quantization bin (recall that $u_1=\mathbb{E}[m_1|x_1]$), and a fully-informative equilibrium exists on the $y$-dimension since there is no bias on that dimension. It is observed that the encoder policy is linear due to the unbiased property of the $y$-dimension.
\begin{figure} [ht]
	\begin{center}
		\includegraphics[scale=0.5]{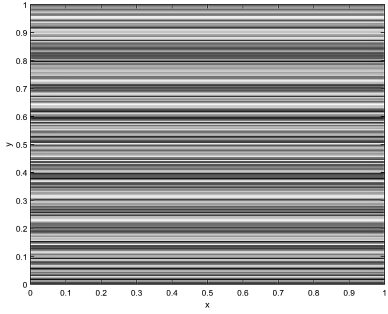}
	\end{center}
	\caption{There is 1 quantization level on the $x$-dimension and 200 quantization levels on the $y$-dimension. The number of quantization levels on the $y$-dimension can be arbitrarily chosen (since $\vb$ is orthogonal to that dimension). As the number of levels goes to infinity, this construction converges to the structure of a linear equilibrium.}
	\label{figure:equilibriumInf2}
\end{figure}
\end{IEEEproof}

\begin{figure} [ht]
	\begin{center}
		\includegraphics[scale=0.5]{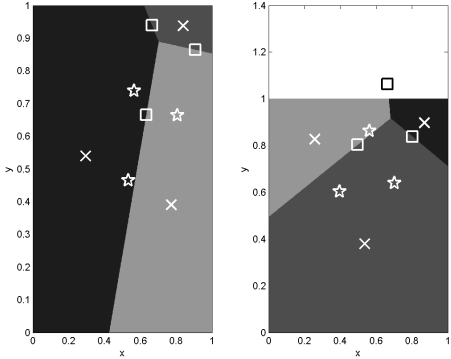}
	\end{center}
	\caption{Sample equilibria in 2D with $\vb_x=0.1$ and $\vb_y=0.2$ where the crosses indicate the centroids of the bins, the star indicates the middle point and the square indicates the shifted middle point.}
	\label{figure:equilibriumMultiple}
\end{figure}

Besides linear equilibria, there may be multiple (hence, non-unique) quantized equilibria with finite regions in the multi-dimensional case as illustrated in Fig.~\ref{figure:equilibriumMultiple}.

From the discussion above, it can be deduced that if $\vb$ is orthogonal to the basis vectors or satisfies certain symmetry conditions, then non-discrete or linear equilibria exist. This approach applies also to the $n$-dimensional setup for any $n \in \mathbb{N}$. For example, if the bias vector involves only one nonzero coordinate component and if the source distribution is uniform over an $n$-dimensional unit cube, then full information revelation in all the other coordinates will lead to a non-discrete equilibrium. In particular, if nonzero component of the bias is greater than $0.25$, then there is only one bin in that coordinate and the full information is sent in other coordinates. Furthermore, if the encoder only sends the $0$ variable for the value of the only bin in the coordinate for which the bias has nonzero component, then what we have is indeed a linear policy. 

\subsection{Multi-Dimensional Cheap Talk: Stackelberg Equilibria}
The Stackelberg equilibria in the multi-dimensional cheap talk can be obtained by extending its scalar case; i.e., it is unique and corresponds to a fully revealing (fully informative) encoder policy as in the scalar case. Thus, Theorem~\ref{thm:cheapStackelberg} holds for the multi-dimensional case as well.

\section{Quadratic Signaling Game: Scalar Case}\label{ScalarSignaling}

\begin{figure} [ht]
	\begin{center}
		\includegraphics[scale=0.5]{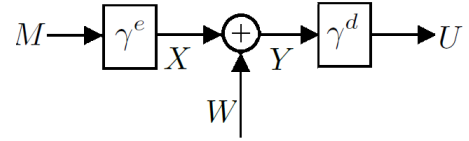}
	\end{center}
	\caption{General system model for noisy case.}
	\label{figure:genScheme}
\end{figure}

The noisy game setup is similar to the noiseless case except that there exists an additive Gaussian noise channel between the encoder and decoder, as shown in Fig.~\ref{figure:genScheme}, and the encoder has a {\it soft} power constraint.

The encoder (DM $1$) encodes a zero-mean Gaussian random variable $M$ and sends the real-valued random variable $X$. During the transmission, the zero mean Gaussian noise with a variance of $\sigma^2$ is added to $X$; hence, the decoder (DM $2$) receives $Y = X + W$. The policy space of DM $1$, $\Gamma^e$, is similarly defined as the policy space in the noiseless case: the set of stochastic kernels from $\mathbb{R}$ to $\mathbb{R}$ (this can be viewed as the measurable subset of the space of all product measures on $\mathbb{R}^2$ with a fixed input marginal, under the weak topology). The policy space of DM $2$, $\Gamma^d$, is the set of stochastic kernels from $\mathbb{R}$ to $\mathbb{R}$. The cost functions of the encoder and the decoder are also slightly modified as follows: DM $1$ aims to minimize
\begin{align*} 
&J^e(\gamma^e,\gamma^d) \\
&\qquad= \int c^e(m,x,u) \gamma^e(\mathrm{d}x|m)\gamma^d(\mathrm{d}u|y) P(\mathrm{d}y|x) P(\mathrm{d}m), 
\end{align*}
whereas DM $2$ aims to minimize
\begin{align*} 
&J^d(\gamma^e,\gamma^d) \\
&\qquad= \int c^d(m, u) \gamma^e(\mathrm{d}x|m)\gamma^d(\mathrm{d}u|y) P(\mathrm{d}y|x) P(\mathrm{d}m), 
\end{align*}
where $P(\mathrm{d}y|x) = P(W \in \mathrm{d}y-x)$ with $W \sim \mathcal{N}(0,\sigma^2)$. The cost functions are modified as $c^e\left(m,x,u\right) = \left(m-u-b\right)^2 + \lambda x^2$ and $c^d\left(m,u\right) = \left(m-u\right)^2$.
Note that a power constraint with an associated multiplier is appended to the cost function of the encoder, which corresponds to power limitation for transmitters in practice. If $\lambda =0$, this corresponds to the setup with no power constraint at the encoder. Here, as earlier, the signaling game problem is investigated where the encoder and the decoder are deterministic; i.e., $\gamma^e(\mathrm{d}x|m)=\mathbbm{1}_{\{f^e(m)\in \mathrm{d}x\}}$ and $\gamma^d(\mathrm{d}u|y)=\mathbbm{1}_{\{f^d(y)\in \mathrm{d}u\}}$ where $f^e(m)$ and $f^d(y)$ are some deterministic functions of the encoder and decoder, respectively.

\subsection{A Supporting Result}
Suppose that there is an equilibrium with an arbitrary policy leading to finite (at least two), countably infinite or uncountably infinite equilibrium bins. Let two of these bins be $\mathcal{B}^{\alpha}$ and $\mathcal{B}^{\beta}$. Also let $m^{\alpha}$ indicate any point in $\mathcal{B}^{\alpha}$; i.e., $m^{\alpha} \in \mathcal{B}^{\alpha}$; and the encoder encodes $m^{\alpha}$ to $x^{\alpha}$ and sends to the decoder. Similarly, let $m^{\beta}$ represent any point in $\mathcal{B}^{\beta}$; i.e., $m^{\beta} \in \mathcal{B}^{\beta}$; and the encoder encodes $m^{\beta}$ to $x^{\beta}$ and sends to the decoder. Without any loss of generality, we can assume that $m^{\alpha}<m^{\beta}$. The decoder chooses the action $u=\mathbb{E}\left[m\middle\vert y\right]$ (MMSE rule). Let $F(m,x)$ be the encoder cost when message $m$ is encoded as $x$; i.e.,
\begin{align*}
F(m,x) &= \int_{y}^{} \! p\left(\gamma^d\left(y\right)=u\middle\vert\gamma^e\left(m\right)=x\right)  \nn\\
& \quad\quad\quad\quad\quad\quad \times\Big((m-u-b)^2+\lambda x^2\Big)  \mathrm{d}y \,.
\end{align*}  
Then the equilibrium definitions from the view of the encoder requires $F(m^\alpha,x^\alpha) \leq F(m^\alpha,x^\beta)$ and $F(m^\beta,x^\beta) \leq F(m^\beta,x^\alpha)$. Now let $G(m)=F(m,x^\alpha)-F(m,x^\beta)$. If it can be shown that $G(m)$ is a continuous function of $m$ on the interval $[m^{\alpha},m^{\beta}]$, then it can be deduced that $\exists \, \overline{m} \in [m^{\alpha},m^{\beta}]$ such that $G(\overline{m})=0$ by the Mean Value Theorem since $G(m^\alpha) \leq 0$ and $G(m^\beta) \geq 0$.

\begin{proposition}
$G(m)$ is a continuous function of $m$ on the interval $[m^{\alpha},m^{\beta}]$.
\label{claim:contOfG}
\end{proposition}
\begin{IEEEproof}
It suffices to show that $F(m,x)$ is continuous in $m$. Let $\{m_n\}$ be a sequence which converges to $m$. Recall that $(m_n-u-b)^2 \leq 2 m_n^2 + 2 (u+b)^2 < \infty$ since $m$ is bounded from above and below ($m \in [m^{\alpha},m^{\beta}]$), $b$ is a finite bias and $\mathbb{E}[u^2]=\mathbb{E}[(\gamma^d\left(y\right))^2]<\infty$ (note that any finite cost $\mathbb{E}[(m-u^2)]$ inevitably leads to a finite $\mathbb{E}[u^2]$ since $\mathbb{E}[u^2] = \mathbb{E}[(m+u-m)^2] \leq 2 \mathbb{E}[m^2] + 2\mathbb{E}[(m-u^2)]<\infty$). Then, by the dominated convergence theorem, 
\begin{align*}
&\lim_{n\to\infty} F(m_n,x) =  \lim_{n\to\infty} \int_{y}^{} \! p\left(\gamma^d\left(y\right)=u\middle\vert\gamma^e\left(m_n\right)=x\right)  \nn\\
& \quad\quad\quad\quad\quad\quad\quad \times\Big((m_n-u-b)^2+\lambda x^2\Big)  \mathrm{d}y \nn\\
&\qquad\qquad\quad\quad\;\;= \int_{y}^{} \! p\left(\gamma^d\left(y\right)=u\middle\vert\gamma^e\left(m\right)=x\right)  \nn\\
& \quad\quad\quad\quad\quad\quad\quad \times\Big((m-u-b)^2+\lambda x^2\Big)  \mathrm{d}y = F(m,x)\,,
\end{align*}
which shows the continuity of $F(\cdot,x)$ in the interval $(m^{\alpha},m^{\beta})$.
\end{IEEEproof}
\par From Proposition~\ref{claim:contOfG}, $\exists \, \overline{m} \in [m^{\alpha},m^{\beta}]$ such that $G(\overline{m})=0$ which implies $ F(\overline{m},x^\alpha)=F(\overline{m},x^\beta)$. Then 
\begin{eqnarray}
&\int_{y}^{} \! p\left(\gamma^d\left(y\right)=u\middle\vert\gamma^e\left(\overline{m}\right)=x^\alpha\right)   \nonumber\\
& \quad \qquad \qquad \qquad  \qquad  \times((\overline{m}-u-b)^2+\lambda (x^\alpha)^2)  \mathrm{d}y \nonumber\\
&= \int_{y}^{} \! p\left(\gamma^d\left(y\right)=u\middle\vert\gamma^e\left(\overline{m}\right)=x^\beta\right) \nonumber\\
& \quad \qquad \qquad \qquad  \qquad  \times((\overline{m}-u-b)^2+\lambda (x^\beta)^2)  \mathrm{d}y.\nonumber
\end{eqnarray}
As a result,
\begin{align}
\overline{m} &={\mathbb{E}[(\gamma^d(y))^2|x^\beta] - \mathbb{E}[(\gamma^d(y))^2|x^\alpha] \over  2\left(\mathbb{E}[\gamma^d(y)|x^\beta] - \mathbb{E}[\gamma^d(y)|x^\alpha]\right)}\nn\\ 
&\qquad\qquad+{\lambda\left((x^\beta)^2-(x^\alpha)^2\right)  \over  2\left(\mathbb{E}[\gamma^d(y)|x^\beta] - \mathbb{E}[\gamma^d(y)|x^\alpha]\right)}+ b
\label{eq:noisyEquilibriumEq}
\end{align}
is obtained. Recall that the arguments in Theorem~\ref{noiselessDiscrete} cannot be applied here because of the presence of noise. However, when there is noise in a communication channel, the relation between $\mathbb{E}[u|x]$, $\mathbb{E}[u^2|x]$ and $\overline{m}$ can be constructed as in \eqref{eq:noisyEquilibriumEq}. 

\subsection{Existence and Uniqueness of Informative Equilibria and Affine Equilibria}
We first note that Proposition~\ref{prop:nonInEq} is valid also in the noisy formulation; i.e. a non-informative (babbling) equilibrium is an equilibrium for the noisy signaling game, since the appended power constraint is always positive. The following holds:
\begin{thm} \label{thm:noisyScalarLambda}
	\begin{enumerate}
		\item Let $0<\lambda < {\mathbb{E}[m^2] \over \mathbb{E}[w^2]}$. For any $b \in \mathbb{R}$, there exists a unique informative affine equilibrium.
		\item If $\lambda \geq  {\mathbb{E}[m^2] \over \mathbb{E}[w^2]}$, there does not exist an informative (affine or non-linear or even randomized) equilibrium. The only equilibrium is the non-informative one.
		\item If $\lambda = 0$, there exists no informative equilibrium with affine policies.
	\end{enumerate}  
\end{thm}
Before presenting the proof, we make the following remark.
\begin{remark}
The expression ${\mathbb{E}[m^2] \over \mathbb{E}[w^2]}$ defines a quantity which determines the Shannon-theoretic capacity of the channel given a signal energy constraint at the encoder. This can be interpreted as Signal-to-Noise Ratio (SNR) of the received signal, which is related to the channel attenuation coefficient. If the multiplier of the signal $\lambda$ in the cost function is greater than $\frac{\mathbb{E}[m^2]}{\mathbb{E}[w^2]}$, it will not be rational for the encoder to send any signal at all under any equilibrium. 
\end{remark}
\begin{IEEEproof}
	\begin{enumerate}
		\item If the encoder is linear (affine), the decoder, as an MMSE decoder for a Gaussian source over a Gaussian channel, is linear (affine); this follows from the property of the conditional expectation for jointly Gaussian random variables. Suppose on the other hand that the decoder is affine so that $u = \gamma^d(y) = Ky + L$ and the encoder policy is $x = \gamma^e(m)$. We will show that the encoder is also affine in this case: With $y = \gamma^e(m)+w$, it follows that $u=K\gamma^e(m)+Kw+L$. By completing the square, the optimal cost of the encoder can be written as
		\begin{align*}
		& J^{*,e} = \min_{x = \gamma^e(m)} \mathbb{E} [(m-u-b)^2+\lambda x^2] \nn\\
		& = \min_{\gamma^e(m)} (K^2+\lambda) \mathbb{E}\Big[\Big(\gamma^e(m)-\frac{(m-L-b)K}{K^2+\lambda}\Big)^2\Big] \nn\\
		&\quad +\frac{\lambda}{K^2+\lambda}\Big(\mathbb{E}[m^2]+(L+b)^2\Big) +K^2\mathbb{E}[w^2]\,.
		\end{align*}
		Hence, the optimal $\gamma^e(m)$ can be chosen as 
		\begin{align}
		\gamma^{*,e}(m)=\frac{(m-L-b)K}{K^2+\lambda} = \frac{(m-L-b)}{K+\lambda/K}\,,
		\label{eq:encoderLinearEquation}
		\end{align}
		and the minimum encoder cost is obtained as 
		\begin{align}
		J^{*,e} = \frac{\lambda}{K^2+\lambda}\Big(\mathbb{E}[m^2]+(L+b)^2\Big) +K^2\mathbb{E}[w^2]\,.
		\label{eq:encoderCost}
		\end{align}
		Recall that \eqref{eq:encoderLinearEquation} implies that an optimal encoder policy for a Gaussian source over a Gaussian channel is an affine policy if the decoder policy is chosen as affine. 		
		We now wish to see if these sets of policies satisfy a fixed point equation. If the decoder has an affine policy, it is proved that the optimal policy of the encoder is also affine :
		\begin{align}
		\gamma^e(m) = Am+C = \Big(\frac{1}{K+\lambda/K}\Big)m+\Big(\frac{-L-b}{K+\lambda/K}\Big)\,.
		\label{eq:linearPolicyEq}
		\end{align}
		On the other hand, with the given affine encoding policy $x =\gamma^e(m) = Am+C$, the optimal decoder policy would be 
		\begin{align}
		\gamma^d(y) = Ky+L = {A \mathbb{E}[m^2] \over A^2 \mathbb{E}[m^2] + \mathbb{E}[w^2]} (y-C)\,.
		\label{eq:decoderPolicyEq}
		\end{align}
		By combining these, we obtain $(K^2+\lambda)^2 \mathbb{E}[w^2] = \lambda \mathbb{E}[m^2]$ by assuming $A\neq0$; which implies $K^2=\sqrt{\frac{\lambda \mathbb{E}[m^2]}{\mathbb{E}[w^2]}}-\lambda$. If we combine the equations above by using $A$, and  define the resulting mapping as $T(A)$, we obtain
		\begin{align}
		A &= {{A \over A^2 + \mathbb{E}[w^2]/\mathbb{E}[m^2]} \over \left({A \over A^2 + \mathbb{E}[w^2]/\mathbb{E}[m^2]}\right)^2+\lambda} \triangleq T(A)\,.
		\label{eq:fixedPointEq}
		\end{align}
		Note now that
		\begin{align*}
		A\geq1 \Rightarrow& {A \over A^2 + {\mathbb{E}[w^2]\over\mathbb{E}[m^2]}}<1 
		\Rightarrow T(A) < {1 \over \lambda}\,, \nn\\
		A<1 \Rightarrow& {A \over A^2 + {\mathbb{E}[w^2]\over\mathbb{E}[m^2]}}<{\mathbb{E}[m^2] \over \mathbb{E}[w^2]} 
		\Rightarrow T(A) 
		< {{\mathbb{E}[m^2]\over\mathbb{E}[w^2]} \over \lambda}\,,
		\end{align*}
		which implies that the mapping defined by $T(A)=A$ can be viewed as a continuous function mapping the compact convex set $[0, {\max(\mathbb{E}[m^2]/\mathbb{E}[w^2], 1) / \lambda}]$ to itself. Therefore, by Brouwer's fixed point theorem \cite{InfiniteDimensionalAnalysis}, there exists $A=T(A)$.
		Indeed, we can find nonzero $A$ for every $0<\lambda<{\mathbb{E}[m^2] \over \mathbb{E}[w^2]}$.\footnote{Recall that if $A\neq0$ and $0<\lambda<{\mathbb{E}[m^2] \over \mathbb{E}[w^2]}$, we have $K^2=\sqrt{\frac{\lambda \mathbb{E}[m^2]}{\mathbb{E}[w^2]}}-\lambda$, which implies $A=\frac{1}{K+\lambda/K}=\pm\sqrt{\sqrt{{\mathbb{E}[w^2] \over \lambda \mathbb{E}[m^2]}}-{\mathbb{E}[w^2] \over \mathbb{E}[m^2]}}$.} After finding $A$, the values for $K$, $C$ and $L$ can also be obtained based on the equilibrium equations in \eqref{eq:linearPolicyEq} and \eqref{eq:decoderPolicyEq}.
		For the uniqueness of an informative fixed point, suppose that there are two different nonzero fixed points: $A_1=T(A_1)$ and $A_2=T(A_2)$ and let $\gamma=\mathbb{E}[w^2]/\mathbb{E}[m^2]$ for simplicity. Then ${A_1 / T(A_1)} = {A_2 / T(A_2) }$ implies
		\begin{align*}
		& {A_1^2 \over A_1^2+\gamma} + \lambda(A_1^2+\gamma) = {A_2^2 \over A_2^2+\gamma} + \lambda(A_2^2+\gamma) \nn\\
		& \Rightarrow (A_1^2-A_2^2)\Big({\gamma \over (A_1^2+\gamma)(A_2^2+\gamma)}+\lambda\Big)=0 \,.
		\end{align*}
		Hence, $|A_1|=|A_2|$ is obtained, and since the mapping is defined from $[0, {\max(\mathbb{E}[m^2]/\mathbb{E}[w^2], 1) / \lambda}]$ to itself, the nonzero fixed point is unique. Then the encoder may choose the nonzero fixed point for the informative equilibirum if it results in a lower cost than the non-informative equilibrium (due to the cost of communication, an informative equilibrium is not always beneficial to the encoder compared to the non-informative one).
		\item Let $\lambda \geq \mathbb{E}[m^2] / \mathbb{E}[w^2]$ and suppose that we are in an equilibrium. Then, the encoder cost $J^e = \mathbb{E}[(m-u-b)^2+\lambda x^2]$ reduces to $J^e = \mathbb{E}[(m-u)^2]+\lambda \mathbb{E}[x^2]+b^2$ since the decoder in an equilibrium always chooses $u=\mathbb{E}[m|y]$. Through $P=\mathbb{E}[x^2]$, the following analysis leads to a lower bound on the encoder cost:
		\begin{align}
		J^e &= b^2 + \lambda \mathbb{E}[x^2] + \mathbb{E}[(m-u)^2] \nonumber \\
		& \overset{(a)}{\geq} b^2 + \lambda P + \mathbb{E}[m^2]{\rm{e}}^{-2 \sup I(X;Y)} \nn\\
		&  = b^2 + \lambda P + \mathbb{E}[m^2]{\rm{e}}^{-2 \frac{1}{2} \log\big(1+\frac{P}{\mathbb{E}[w^2]}\big)} \nonumber \\
		& = b^2 + \lambda P + \frac{\mathbb{E}[m^2]}{1+P /\mathbb{E}[w^2]}\,.		\label{eq:itInequality}
		\end{align}
Here, (a) follows from a rate-distortion theoretic bound through the data-processing inequality (see for example p. 96 of \cite{YukselBasarBook}). However, it follows that when $\lambda \geq \mathbb{E}[m^2] / \mathbb{E}[w^2]$, \eqref{eq:itInequality} is minimized at $P=0$; that is, the encoder does not signal any output. Hence, the encoder engages in a non-informative equilibrium and the minimum cost becomes $\mathbb{E}[m^2]+b^2$ at this non-informative equilibrium.
		\item It is proved that an optimal encoder is affine such that $x =\gamma^e(m) = Am+C$ when the decoder is affine, that is, $u = \gamma^d(y) = Ky + L$. Then, by inserting $\lambda=0$ to \eqref{eq:noisyEquilibriumEq}, $\overline{m}$ is obtained as $\overline{m} = KA{\frac{(m^\alpha+m^\beta)}{2}}+KC+L+b$.
		This holds for all $m^\alpha$ and $m^\beta$ with $m^{\alpha} \leq \overline{m} \leq m^{\beta}$. Thus, if the distance between $m^\alpha$ and $m^\beta$ is made arbitrarily small, then it must be that $KA=1$ and $KC+L+b=0$. On the other hand, it was shown that an optimal decoder policy is affine if an encoder is affine in \eqref{eq:decoderPolicyEq}. By combining $KA=1$ and $K = {A \mathbb{E}[m^2] \over A^2 \mathbb{E}[m^2] + \mathbb{E}[w^2]}$, it follows that a real-valued solution does not exist for any given affine coding parameter. 
	\end{enumerate}
\end{IEEEproof}

\begin{remark}
	Note that, from \eqref{eq:linearPolicyEq} and \eqref{eq:decoderPolicyEq}, we have $A={1 \over K+\lambda/K}$, $K={A \mathbb{E}[m^2] \over A^2 \mathbb{E}[m^2] + \mathbb{E}[w^2]}$, $L=-KC$ and $Ab=(AK-1)C$. From these equalities, we observe the following:
\begin{enumerate}
	\item when $\lambda=0$, it is shown in Theorem~\ref{thm:noisyScalarLambda} that there is not any fixed point solution to \eqref{eq:fixedPointEq}. However, if there is not a noisy channel between the encoder and the decoder; i.e., the noise variance is zero ($\mathbb{E}[w^2]=0$), then \eqref{eq:fixedPointEq} has a fixed point solution. Even when \eqref{eq:fixedPointEq} has a fixed point solution $A$, \eqref{eq:linearPolicyEq} and \eqref{eq:decoderPolicyEq} cannot hold together unless $b=0$.
	\item when the noise variance is zero ($\mathbb{E}[w^2]=0$), there is not any fixed point solution to \eqref{eq:fixedPointEq} unless $\lambda=0$. Even when \eqref{eq:fixedPointEq} has a fixed point solution $A$,  \eqref{eq:linearPolicyEq} and \eqref{eq:decoderPolicyEq} cannot hold together unless $b=0$.
	\item when $\lambda = 0$ and the noise variance is zero ($\mathbb{E}[w^2]=0$); the consistency of \eqref{eq:linearPolicyEq} and \eqref{eq:decoderPolicyEq} can be satisfied if only if $b=0$. Hence, if $b \neq 0$, there cannot be a affine equilibrium; the equilibrium has to be discrete due to Theorem~\ref{noiselessDiscrete}. 
\end{enumerate} 
Thus, if either $\lambda$ or $\mathbb{E}[w^2]$ is $0$, an affine equilibrium exists only if $\lambda$, $\mathbb{E}[w^2]$ and $b$ are all $0$.
\end{remark}

\subsection{Price of Anarchy and Comparison with Socially Optimal Cost}

In a game theoretic setup, the encoder and the decoder try to minimize their individual costs, thus the game theoretic cost can be found as $\min_{\gamma^e}J^e + \min_{\gamma^d}J^d$. If the encoder and the decoder work together to minimize the total cost, then the problem can be regarded as a team problem and the resulting cost is a socially optimal cost, which is $\min_{\gamma^e, \gamma^d}(J^e + J^d)$. In the game theoretic setup, because of the selfish behavior of the players, there is some loss from the socially optimal cost, and this loss is measured by the ratio between the game theoretic cost and the socially optimal cost, which was proposed as a price of anarchy \cite{worstCaseEquilibria}. In this part, it will be shown that the game theoretic cost is higher than the socially optimal cost as expected, and the information theoretic lower bounds on the costs and their achievability will be discussed. 
\begin{thm}\label{thm:gameTeamCost}
\begin{enumerate}
\item Let $g_i$ and $g_u$ represent the informative and the non-informative equilibrium game costs, respectively. Then, $g_i=3\sqrt{\lambda \mathbb{E}[m^2] \mathbb{E}[w^2]}+b^2\sqrt{{\mathbb{E}[m^2] \over \lambda \mathbb{E}[w^2]}}-\lambda \mathbb{E}[w^2]$ and $g_u=2\mathbb{E}[m^2]+b^2$. Further, the total cost in the game equilibrium is the following
\begin{align*}
		J^{*,g} = \begin{cases}
			\min \{g_i,\; g_u\}  & \lambda < \mathbb{E}[m^2]/\mathbb{E}[w^2] \\
			g_u & \lambda\geq \mathbb{E}[m^2]/\mathbb{E}[w^2] 
		\end{cases} \,.
\end{align*}
\item Let $t_i$ and $t_u$ represent the informative and the non-informative team costs, respectively. Then, $t_i=2\sqrt{2\lambda \mathbb{E}[m^2]\mathbb{E}[w^2]}+{b^2\over2}-\lambda \mathbb{E}[w^2]$ and $t_u=2\mathbb{E}[m^2]+{b^2\over2}$. Further, the socially optimal cost (the total cost in the team setup) is the following
\begin{align*}
		J^{*,t} = \begin{cases}
			\min \{t_i,\; t_u\} & \lambda < 2 \mathbb{E}[m^2]/\mathbb{E}[w^2] \\
			t_u& \lambda\geq 2 \mathbb{E}[m^2]/\mathbb{E}[w^2] 
		\end{cases}\,.
\end{align*}
\end{enumerate}
\end{thm}
\begin{IEEEproof}
\begin{enumerate}
\item Note from \eqref{eq:linearPolicyEq} and \eqref{eq:decoderPolicyEq} that we have $A={1 \over K+\lambda/K}$, $K={A \mathbb{E}[m^2] \over A^2 \mathbb{E}[m^2] + \mathbb{E}[w^2]}$, $L=-KC$ and $Ab=C(AK-1)$. Also we have  $(K^2+\lambda)^2 \mathbb{E}[w^2] = \lambda \mathbb{E}[m^2]$ which implies $K^2=\sqrt{\frac{\lambda \mathbb{E}[m^2]}{\mathbb{E}[w^2]}}-\lambda$ and $\lambda<\mathbb{E}[m^2]/\mathbb{E}[w^2]$ for nonzero $A$. Recall that if $\lambda\geq \mathbb{E}[m^2]/\mathbb{E}[w^2]$, then $A=C=K=L=0$, which implies the non-existence of the informative linear  (also affine) equilibrium. Thus, for $\lambda<{\mathbb{E}[m^2] \over \mathbb{E}[w^2]}$, by using $K^2=\sqrt{\frac{\lambda \mathbb{E}[m^2]}{\mathbb{E}[w^2]}}-\lambda$, $A={1 \over K+\lambda/K}$, $C = {Ab \over AK-1}$ and $L+b=-{C\over A}$ in \eqref{eq:encoderCost}, we have
\[J^{*,e} = 2\sqrt{\lambda \mathbb{E}[m^2] \mathbb{E}[w^2]}+b^2\sqrt{{\mathbb{E}[m^2] \over \lambda \mathbb{E}[w^2]}}-\lambda \mathbb{E}[w^2]\,. \]
Now recall that the optimal decoder policy is $u^* = \mathbb{E}[m|(y=Am+C+w)] = {A \mathbb{E}[m^2] \over A^2 \mathbb{E}[m^2] + \mathbb{E}[w^2]} (y-C)$, and we have $\sigma_e^2=\sigma_x^2-{\sigma_{xy}^2\over\sigma_y^2}$ where $e=x-\mathbb{E}[x|y]$. In this case, $x\rightarrow m$, $y\rightarrow y$, $\sigma_x^2 \rightarrow \mathbb{E}[m^2]$, $\sigma_{xy}\rightarrow A\mathbb{E}[m^2]$ and $\sigma_y^2\rightarrow A^2\mathbb{E}[m^2]+\mathbb{E}[w^2]$. Thus, we have
\begin{align*}
J^{*,d} &= \min_{u = \gamma^d(y)} \mathbb{E} [(m-u)^2] = \mathbb{E} [(m-\mathbb{E}[m|y])^2] \nn\\
&=\sigma_m^2-{\sigma_{my}^2\over\sigma_y^2}=\mathbb{E}[m^2]-{A^2(\mathbb{E}[m^2])^2\over A^2\mathbb{E}[m^2]+\mathbb{E}[w^2]} \nn\\
&= \sqrt{\lambda \mathbb{E}[m^2] \mathbb{E}[w^2]}\,.
\end{align*}
As a result, the game theoretic cost at the equilibrium is found as 
\begin{align}
J^{*,g}=3\sqrt{\lambda \mathbb{E}[m^2] \mathbb{E}[w^2]}+b^2\sqrt{{\mathbb{E}[m^2] \over \lambda \mathbb{E}[w^2]}}-\lambda \mathbb{E}[w^2]\,.
\label{eq:gameCost}
\end{align}
Recall that, if $\lambda\geq \mathbb{E}[m^2]/\mathbb{E}[w^2]$, then $J^{*,e}=\mathbb{E}[m^2]+b^2$ and $J^{*,d}=\mathbb{E}[m^2]$; hence, $J^{*,g}=2\mathbb{E}[m^2]+b^2$. If there were no cost of communication (consider the cheap talk; i.e., remove $\lambda x^2$ from the encoder cost function), then one could say that the informative equilibria would always be beneficial to both the encoder and the decoder; however, due to the cost of communication, an informative equilibrium is not always beneficial to the encoder when compared with the non-informative one (i.e., for $\lambda< \mathbb{E}[m^2]/\mathbb{E}[w^2]$, it does not always hold that $2\sqrt{\lambda \mathbb{E}[m^2] \mathbb{E}[w^2]}+b^2\sqrt{{\mathbb{E}[m^2] \over \lambda \mathbb{E}[w^2]}}-\lambda \mathbb{E}[w^2] < \mathbb{E}[m^2]+b^2$). For the receiver, however, information never hurts the performance and the informative equilibria are more desirable (i.e., for $\lambda< \mathbb{E}[m^2]/\mathbb{E}[w^2]$, the inequality $\sqrt{\lambda \mathbb{E}[m^2] \mathbb{E}[w^2]}<\mathbb{E}[m^2]$ always holds). As a result, one can expect a non-informative equilibrium even if $\lambda< \mathbb{E}[m^2]/\mathbb{E}[w^2]$.
\item The part below aims to construct the socially optimal affine setup. In this part, $J^{e,t}$ represents the team cost minimized over the encoder policies for a given decoder policy, $J^{d,t}$ represents the team cost minimized over the decoder policies for a given encoder policy, and $J^{*,t}$ represents the optimum team cost; i.e., minimization over all affine encoding and decoding policies as follows:
\[J^{*,t} = \min_{x = \gamma^e(m), u = \gamma^d(y)} \mathbb{E} [(m-u-b)^2+\lambda x^2+(m-u)^2]\,. \]
\par Similar to the game theoretic analysis above, with the given affine encoding policy $x =\gamma^e(m) = Am+C$ (then $y=x+w=Am+C+w$), the optimal decoder policy can be found as follows (by completing the square):
\begin{align*}
J^{d,t} &= \min_{u = \gamma^d(y)} \mathbb{E} [(m-u-b)^2+\lambda x^2+(m-u)^2] \nn\\
&= \min_{u = \gamma^d(y)} 2 \mathbb{E}\Big[ (m-u - {b\over2})^2 + {b^2\over4} + \lambda {x^2\over2}\Big]\,.
\end{align*}
Hence the optimal decoder policy can be chosen as $\gamma^{d,t}(y)=\mathbb{E}[m-{b\over2}\,|\,y]$. Due to the joint Gaussanity of $m$ and $y$, the minimizer decoder policy is affine:
\begin{align}
\gamma^{d,t}(y) = Ky+L = {A \mathbb{E}[m^2] \over A^2 \mathbb{E}[m^2] + \mathbb{E}[w^2]} (y-C) - {b\over2}\,. \label{eq:teamDecoder}
\end{align}
\par Similar to the game theoretic analysis above, for any affine decoder policy $\gamma^d(y) = Ky+L$ with $y=\gamma^e(m)+w$, the optimal encoder policy for the team setup can be obtained as follows (by completing the square):
\begin{align*}
& J^{e,t} = \min_{x = \gamma^e(m)} \mathbb{E} [(m-u-b)^2+\lambda x^2+(m-u)^2] \nn\\
& = \min_{\gamma^e(m)} (2K^2+\lambda) \mathbb{E}\Big[\Big(\gamma^e(m)-\frac{(2m-2L-b)K}{2K^2+\lambda}\Big)^2\Big] \nn\\
& \;+\frac{b^2K^2 + \lambda\left(2\mathbb{E}[m^2]+(L+b)^2+L^2\right)}{2K^2+\lambda}+2K^2\mathbb{E}[w^2]\,.\nn
\end{align*} 
Hence, the optimal encoder $\gamma^e(m)$ is
\begin{align}
\gamma^{e,t}(m)=Am+C=\frac{(2m-2L-b)}{2K+\lambda/K}\,,
\label{eq:teamEncoder}
\end{align}
and the minimum team cost is obtained as
\begin{align}
J^{*,t} &= \frac{b^2K^2 + \lambda\left(2\mathbb{E}[m^2]+(L+b)^2+L^2\right)}{2K^2+\lambda}\nn\\
&\qquad\qquad\qquad\qquad\qquad+2K^2\mathbb{E}[w^2]\,.
\label{eq:teamCost}
\end{align}
This implies that, in the team setup, an optimal encoder policy for a Gaussian source over a Gaussian channel is a affine policy if the decoder policy is chosen as affine. 

In order to achieve the socially optimal cost $J^{*,t}$, the optimal encoder policy $\gamma^{e^*,t}(m)$ and the optimal decoder policy $\gamma^{d^*,t}(y)$ must satify the following equalities by  \eqref{eq:teamDecoder} and \eqref{eq:teamEncoder}:
\begin{align*}
A &= {2\over2K+\lambda/K}\,,\quad \quad K = {A \mathbb{E}[m^2] \over A^2 \mathbb{E}[m^2] + \mathbb{E}[w^2]}\,,\nn\\
C &= {A\over2}(-2L-b) = -AL-{Ab\over2}\,,\quad L = -KC-{b\over2} \nn\\
\Rightarrow & C=-A\left(-KC-{b\over2}\right)-{Ab\over2}=AKC\,.
\end{align*}
Here, either $AK=1$ or $C=0$. If $AK=1$, then $\mathbb{E}[w^2] = 0$ which contradicts with the noise assumption. Then $C=0$ and $L=-b/2$. By using the equalities for $A$ and $K$ above, one can obtain $2(K^2+\lambda/2)^2 \mathbb{E}[w^2] = \lambda \mathbb{E}[m^2]$ by assuming $A\neq0$; which implies $K^2=\sqrt{\frac{\lambda \mathbb{E}[m^2]}{2\mathbb{E}[w^2]}}-\frac{\lambda}{2}$. Since $K^2$ is positive, $\lambda$ cannot be greater than $2\mathbb{E}[m^2] \over \mathbb{E}[w^2]$; otherwise, because of our assumption, $A$ must be equal to $0$ which implies that $K=0$, and there does not exist an informative affine team setup. Then $K^2=\sqrt{\frac{\lambda \mathbb{E}[m^2]}{2\mathbb{E}[w^2]}}-\frac{\lambda}{2}$ and $\lambda<2\mathbb{E}[m^2]/\mathbb{E}[w^2]$ for nonzero $A$. Thus, for $\lambda<{2\mathbb{E}[m^2] \over \mathbb{E}[w^2]}$, by using $K^2=\sqrt{\frac{\lambda \mathbb{E}[m^2]}{2\mathbb{E}[w^2]}}-\frac{\lambda}{2}$, $A={2K \over 2K^2+\lambda}$, $C = 0$ and $L=-{b\over2}$ in \eqref{eq:teamCost}, we have
\begin{align}
J^{*,t} = 2\sqrt{2\lambda \mathbb{E}[m^2]\mathbb{E}[w^2]}+{b^2\over2}-\lambda \mathbb{E}[w^2] m\,.
\label{eq:teamCostLast}
\end{align}
Recall that, if $\lambda\geq 2\mathbb{E}[m^2]/\mathbb{E}[w^2]$, then $J^{*,t}=2\mathbb{E}[m^2]+{b^2\over2}$. Similar to the game theoretic setup, due to the cost of the communication, the encoder and the decoder may prefer the non-informative equilibrium over the informative one (if $2\sqrt{2\lambda \mathbb{E}[m^2]\mathbb{E}[w^2]}+{b^2\over2}-\lambda \mathbb{E}[w^2]>2\mathbb{E}[m^2]+{b^2\over2}$).
\end{enumerate}
\end{IEEEproof}
\begin{thm}
The price of anarchy is always larger than 1, i.e., the sum of the costs under any Nash equilibria is always larger than the socially optimal cost.
\label{thm:priceAnarchy}
\end{thm}
\begin{IEEEproof}
By Theorem~\ref{thm:gameTeamCost}, we have the following
\begin{align*}
J^{*,g} = \begin{cases}
\min \{g_i,\; g_u\}  & \lambda < \mathbb{E}[m^2]/\mathbb{E}[w^2] \\
g_u & \lambda\geq \mathbb{E}[m^2]/\mathbb{E}[w^2] 
\end{cases}\;,
\\
J^{*,t} = \begin{cases}
\min \{t_i,\; t_u\} & \lambda < 2 \mathbb{E}[m^2]/\mathbb{E}[w^2] \\
t_u& \lambda\geq 2 \mathbb{E}[m^2]/\mathbb{E}[w^2] 
\end{cases}\;.
\end{align*} 
Notice that we have $t_i<g_i$ for $\lambda < \mathbb{E}[m^2]/\mathbb{E}[w^2]$ and $t_u<g_u$ always. Consider the following cases:
\begin{enumerate}
\item \underline{$0<\lambda < \mathbb{E}[m^2]/\mathbb{E}[w^2]$ :} There are four cases to be considered:
\begin{enumerate}
\item \underline{$\min \{g_i,\; g_u\}=g_i$ and $\min \{t_i,\; t_u\}=t_i$:} Since $t_i<g_i$, $J^{*,t}<J^{*,g}$ is satisfied.
\item \underline{$\min \{g_i,\; g_u\}=g_i$ and $\min \{t_i,\; t_u\}=t_u$:} Since $t_u<t_i<g_i$, $J^{*,t}<J^{*,g}$ is satisfied.  
\item \underline{$\min \{g_i,\; g_u\}=g_u$ and $\min \{t_i,\; t_u\}=t_i$:} Since $t_i<t_u<g_u<g_i$,  $J^{*,t}<J^{*,g}$ is satisfied.
\item \underline{$\min \{g_i,\; g_u\}=g_u$ and $\min \{t_i,\; t_u\}=t_u$:} Since $t_u<g_u$, $J^{*,t}<J^{*,g}$ is satisfied.
\end{enumerate} 
\item \underline{$\mathbb{E}[m^2]/\mathbb{E}[w^2]\leq\lambda < 2\mathbb{E}[m^2]/\mathbb{E}[w^2]$ :} There are two cases to be considered:
\begin{enumerate}
\item \underline{$\min \{t_i,\; t_u\}=t_i$:} Since $t_i<t_u<g_u$, $J^{*,t}<J^{*,g}$ is satisfied.
\item \underline{$\min \{t_i,\; t_u\}=t_u$:} Since $t_u<g_u$, $J^{*,t}<J^{*,g}$ is satisfied. 
\end{enumerate}
\item \underline{$\lambda \geq 2\mathbb{E}[m^2]/\mathbb{E}[w^2]$ :} Since $t_u<g_u$, $J^{*,t}<J^{*,g}$ is satisfied. 
\end{enumerate}
Hence, one can observe that $J^{*,g}>J^{*,t}$ always holds, which shows that the price of anarchy is greater than 1, i.e., the game theoretic cost is always larger than the socially optimal cost.
\end{IEEEproof}

In the following, we discuss information theoretic lower bounds on the performance of equilibria and socially optimal strategies.
\begin{thm}
	\begin{enumerate}
		\item For the game setup, if $\lambda\geq{\mathbb{E}[m^2] \over \mathbb{E}[w^2]}$ (i.e., non-informative equilibria), the information theoretic lower bounds on the costs are achievable.
		\item For the game setup, if $\lambda<{\mathbb{E}[m^2] \over \mathbb{E}[w^2]}$ and $b=0$, then the information theoretic lower bounds on the costs are achievable by linear policies.
		\item For the game setup, if $\lambda<{\mathbb{E}[m^2] \over \mathbb{E}[w^2]}$ and $b\neq0$, the information theoretic lower bounds on the costs are not achievable by affine policies. 
		\item For the team setup, the information theoretic lower bounds on the costs are always (both in the informative and non-informative equilibria) achievable by affine policies.
	\end{enumerate}
	\label{thm:itLowerBound}
\end{thm}

\begin{IEEEproof}
	\begin{enumerate}
		\item Recall that the encoder cost is $J^e = \mathbb{E}[(m-u-b)^2+\lambda x^2]$ and we know that this reduces to $J^e = \mathbb{E}[(m-u)^2]+\lambda \mathbb{E}[x^2]+b^2$ since the decoder always chooses $u=\mathbb{E}[m|y]$. From \eqref{eq:itInequality}, we have a bound on the encoder cost \mbox{$J^e\geq b^2 + \lambda P + \frac{\mathbb{E}[m^2]}{1+P /\mathbb{E}[w^2]}$} where $P=\mathbb{E}[x^2]$ represents the power. This bound is tight when the encoder and the decoder use linear policies leading to jointly Gaussian random variables. For $\lambda < \mathbb{E}[m^2]/\mathbb{E}[w^2]$, a minimizer of this cost is $P^*=\sqrt{{\mathbb{E}[m^2]\mathbb{E}[w^2]\over\lambda}}-\mathbb{E}[w^2]$. If we insert this value into \eqref{eq:itInequality}, we have $J^e \geq 2\sqrt{\lambda \mathbb{E}[m^2]\mathbb{E}[w^2]} + b^2 - \lambda \mathbb{E}[w^2]$.
		By the same reasoning above, we also have $
		J^d = \mathbb{E}[(m-u)^2] \geq \frac{\mathbb{E}[m^2]}{1+{P\over \mathbb{E}[w^2]}} \geq \sqrt{\lambda \mathbb{E}[m^2]\mathbb{E}[w^2]}$.
		Hence, the information theoretic lower bound on the game cost $J^g= J^e+J^d$ is found as
		\begin{align}
		J^g \geq 3\sqrt{\lambda \mathbb{E}[m^2]\mathbb{E}[w^2]} + b^2 - \lambda \mathbb{E}[w^2]\,.
		\label{eq:gameIneq}
		\end{align}
		Through an analysis similar to the one in \cite{YukselBasarBook}, one can see that when $\lambda \geq \mathbb{E}[m^2] / \mathbb{E}[w^2]$, \eqref{eq:itInequality} is minimized at $P=0$ (the encoder does not signal any output); thus we obtain a non-informative equilibrium: The encoder and the decoder do not engage in communications; i.e., $A=0$ and $K=0$ is an equilibrium. In this case the encoder may be considered to be linear, but this is a degenerate coding policy. This implies $J^g \geq 2\mathbb{E}[m^2]+b^2$, and remember that $J^{*,g}=2\mathbb{E}[m^2]+b^2$ when $\lambda \geq \mathbb{E}[m^2] / \mathbb{E}[w^2]$, hence the information theoretic lower bound is achievable in the non-informative equilibria.
		\item From \eqref{eq:gameCost} and \eqref{eq:gameIneq}, it can be deduced that when $b=0$, the lower bound of the encoder cost is achievable by linear policies; i.e., $C=0$ and $L=0$. When $b=0$, the problem corresponds to what is known as a {\it  soft-constrained version of the quadratic signaling problem} where we append the constraint to the cost functional (see page 96 of \cite{YukselBasarBook}). 
		\item If $b\neq 0$, then, from \eqref{eq:gameCost} and \eqref{eq:gameIneq}, one can observe that the lower bound becomes unachievable by affine policies since the power constraint related part of the cost function, $\lambda x^2$, contains $b^2$ related parameters (recall $C = {Ab \over AK-1}$). In this case, by modifying the power from $P$ to $P-C^2$ (which must be positive) in the information theoretic inequalities; i.e., $J^e \geq b^2 + \lambda P + \frac{\mathbb{E}[m^2]}{1+(P-C^2) /\mathbb{E}[w^2]}$, then the minimum game cost is obtained as $J^g \geq 3\sqrt{\lambda \mathbb{E}[m^2] \mathbb{E}[w^2]}+b^2\sqrt{{\mathbb{E}[m^2] \over \lambda \mathbb{E}[w^2]}}-\lambda \mathbb{E}[w^2]$ which is the same cost that is achieved by affine policies.
		\item By following a similar approach to \eqref{eq:itInequality} for finding the lower bound on the socially optimal cost, we can obtain:
		\begin{align*}
		J^t &= \mathbb{E}[(m-u-b)^2+\lambda x^2+(m-u)^2]\nn\\
		&= {b^2\over2} + \lambda \mathbb{E}[x^2] + 2\mathbb{E}\left[\left(m-u-{b\over2}\right)^2\right]\nn\\
		&\overset{(a)}{\geq} {b^2\over2} + \lambda P + \frac{2\mathbb{E}[m^2]}{1+P /\mathbb{E}[w^2]}\,.
		\end{align*}
		Here $(a)$ holds since the decoder chooses $u=\mathbb{E}[m-{b\over2}|y]$ and shifting does not affect the differential entropy. Similar to the previous analysis, a minimizer of this cost is  $P^*=\sqrt{{2\mathbb{E}[m^2]\mathbb{E}[w^2]\over\lambda}}-\mathbb{E}[w^2]$ for $\lambda < 2\mathbb{E}[m^2]/\mathbb{E}[w^2]$. If we insert this value into the total cost, we have 
		\begin{align}
		J^t \geq 2\sqrt{2\lambda \mathbb{E}[m^2]\mathbb{E}[w^2]}+{b^2\over2}- \lambda \mathbb{E}[w^2].
		\label{eq:teamIneq}
		\end{align}
		Recall that, if $\lambda\geq 2\mathbb{E}[m^2]/\mathbb{E}[w^2]$, then $P=0$ becomes the minimizer, hence $J^t \geq 2\mathbb{E}[m^2]+{b^2\over2}$ in the non-informative equilibrium. Remember that $J^{*,t}=2\mathbb{E}[m^2]+{b^2\over2}$ in this case, thus the information theoretic lower bound is achievable in the non-informative equilibria. In addition, from \eqref{eq:teamCostLast} and \eqref{eq:teamIneq}, for $\lambda<2\mathbb{E}[m^2]/\mathbb{E}[w^2]$ (which implies the informative equilibria), it can easily be seen that the information theoretic lower bound is achievable by affine policies (actually the encoder policy is linear and the decoder policy is affine). \hfill$\IEEEQEDhere$
	\end{enumerate}
\end{IEEEproof}

We state the following summary.
\begin{enumerate}
	\item If $\lambda < \mathbb{E}[m^2]/\mathbb{E}[w^2]$ and $b=0$, then the information theoretic lower bound on the game cost is achievable by the linear policies.
	\item If $\lambda < \mathbb{E}[m^2]/\mathbb{E}[w^2]$ and $b\neq0$, then the information theoretic lower bounds on the game cost are not achievable by the affine policies; but they become achievable after slight modification on the power parameter in the information theoretic inequality.
	\item The team cost $J^{*,t}$ in the affine equilibrium is always equal to the information theoretic lower bound on the team cost.
	\item The price of anarchy is always greater than 1: The socially optimal cost is always lower than the cost in any equilibrium.
	\item In the game setup, the non-informative equilibrium may be preferred over the informative equilibrium by the encoder due to the cost of the signal $\lambda x^2$.
\end{enumerate}

\subsection{Stackelberg Setup}
If we consider the Stackelberg setup of the signaling game problem studied in this section; i.e. the encoder knows the policy of the decoder, then it can be shown that the only equilibrium is the linear equilibrium.
\begin{thm}
The only equilibrium in the Stackelberg setup of the signaling game is the linear equilibrium.
\end{thm}
\begin{IEEEproof}
In the proof, first we assume the linear encoding policy and show that the information theoretic lower bound is achieved, then we conclude that the encoder policy must be linear.
Let the encoder policy be $x = \gamma^e(m) = Am + C$. Due to the Stackelberg assumption, the encoder knows that the decoder will use $\gamma^d(y) = u = \mathbb{E}[m|y]$ as an optimal decoder policy to minimize the decoder cost, thus $u=\gamma^d(y) = {A \mathbb{E}[m^2] \over A^2 \mathbb{E}[m^2] + \mathbb{E}[w^2]} (y-C)$ where $y=Am+C+w$. Then the goal of the encoder is to minimize the following:
		\begin{align}
		& J^{*,e} = \min_{x = \gamma^e(m)= Am + C}  \mathbb{E} [(m-u-b)^2+\lambda x^2] \nn\\
		& = \min_{A,\;C}\; \mathbb{E} \Big[\left({m \mathbb{E}[w^2] - A \mathbb{E}[m^2] w \over A^2 \mathbb{E}[m^2] + \mathbb{E}[w^2]}-b\right)^2  +\lambda (Am+C)^2\Big] \nn\\
		& = \min_{A,\;C}\; {\mathbb{E}[m^2] (\mathbb{E}[w^2])^2 + A^2 (\mathbb{E}[m^2])^2 \mathbb{E}[w^2] \over (A^2 \mathbb{E}[m^2] + \mathbb{E}[w^2])^2} +b^2\nn\\
		&\qquad\qquad\qquad\qquad\qquad\qquad\qquad  +\lambda A^2 \mathbb{E}[m^2] + \lambda C^2 \nn\\
		& = \min_{A,\;C}\; {\mathbb{E}[m^2] \mathbb{E}[w^2] \over A^2 \mathbb{E}[m^2] + \mathbb{E}[w^2]} +b^2  +\lambda A^2 \mathbb{E}[m^2] + \lambda C^2 \,.
		\label{eq:stackelbergOptimalCost}
		\end{align}
The optimal encoder cost in \eqref{eq:stackelbergOptimalCost} is achieved for $C^*=0$, and $A^*=0$ for $\lambda\geq \mathbb{E}[m^2]/\mathbb{E}[w^2]$ and $A^*=\sqrt{\sqrt{\frac{\mathbb{E}[w^2]}{\lambda \mathbb{E}[m^2]}}-\frac{\mathbb{E}[w^2]}{\mathbb{E}[m^2]}}$ for $\lambda< \mathbb{E}[m^2]/\mathbb{E}[w^2]$. Then the optimal encoder cost is obtained as $J^{*,e}=\mathbb{E}[m^2]+b^2$ for $\lambda\geq \mathbb{E}[m^2]/\mathbb{E}[w^2]$ and $J^{*,e}=2\sqrt{\lambda \mathbb{E}[m^2]\mathbb{E}[w^2]} + b^2 - \lambda \mathbb{E}[w^2]$ for $\lambda< \mathbb{E}[m^2]/\mathbb{E}[w^2]$. Note that these are the information theoretic lower bounds in the proof of the first part of Theorem~\ref{thm:itLowerBound} and these lower bounds are achieved when the encoder and the decoder use linear policies jointly, which is valid for the current case.
\end{IEEEproof}

\section{Quadratic Signaling Game: Multi-Dimensional Gaussian Noisy Case}

The scalar setup considered in Section~\ref{ScalarSignaling} can be extended to the multi-dimensional Gaussian noisy signaling game problem setup as follows. The encoder (DM $1$) encodes an $n$-dimensional zero-mean Gaussian random variable $\vM$ with the covariance matrix $\Sigma_{\vM}$ and sends the real-valued $n$-dimensional random variable $\vX$. During the transmission, the $n$-dimensional zero-mean Gaussian noise with the covariance matrix $\Sigma_{\vW}$ is added to $\vX$ and the decoder (DM $2$) receives $\vY = \vX + \vW$. The policy space of DM $1$, $\Gamma^e$, and the policy space of DM $2$, $\Gamma^d$, are the set of stochastic kernels from $\mathbb{R}^n$ to $\mathbb{R}^n$. The cost functions of the encoder and the decoder are as follows: DM $1$ aims to minimize
\begin{align*} 
&J^e(\gamma^e,\gamma^d) \\
&\qquad= \int c^e(\vm,\vx,\vu) \gamma^e(\mathrm{d}\vx|\vm)\gamma^d(\mathrm{d}\vu|\vy)P(\mathrm{d}\vy|\vx) P(\mathrm{d}\vm), 
\end{align*}
whereas DM $2$ aims to minimize
\begin{align*} 
&J^d(\gamma^e,\gamma^d) \\
&\qquad= \int c^d(\vm, \vu) \gamma^e(\mathrm{d}\vx|\vm)\gamma^d(\mathrm{d}\vu|\vy)  P(\mathrm{d}\vy|\vx) P(\mathrm{d}\vm), 
\end{align*}
where $P(d\vy|\vx) = P(\vW \in \mathrm{d}\vy-\vx)$ with $\vW \sim \mathcal{N}(0,\Sigma_{\vW})$. The cost functions are $c^e\left(\vm,\vx,\vu\right) = \|\vm-\vu-\vb\|^2 + \lambda \|\vx\|^2$ and $c^d\left(\vm,\vu\right) = \|\vm-\vu\|^2$ where the lengths of the vectors are defined in $L_2$ norm and $\vb$ is the bias vector.. Note that we have appended a power constraint and an associated multiplier. If $\lambda=0$, this corresponds to the setup with no power constraint at the encoder.

\subsection{Affine Equilibria}

\begin{thm}
	\begin{enumerate}
		\item If the encoder is linear (affine), the decoder, as an MMSE decoder for a Gaussian source over a Gaussian channel, is linear (affine).
		\item If the decoder is linear (affine), then an optimal encoder policy for a multi-dimensional Gaussian source over a multi-dimensional Gaussian channel is an affine policy.
		\item An equilibrium encoder policy $\gamma^e (\vm) = A\vm + \vC$ satisfies the equation $A=T(A)$ where $T(A)= \left(FF^T+\lambda I\right)^{-1}F$ and $F=\left(A\Sigma_{\vM}A^T + \Sigma_{\vW}\right)^{-1}A\Sigma_{\vM}$. 
		\item There exists at least one equilibrium.
	\end{enumerate}
\end{thm}
\begin{IEEEproof}
	\begin{enumerate}
		\item Let the affine encoding policy be $\vx =\gamma^e(\vm) = A\vm+\vC$ where $A$ is an $n \times n$ matrix and $\vC$ is an $n \times 1$ vector. Then $\vy = \vx + \vw = A\vm+\vC+\vw$. The optimal cost of the decoder, by the law of the iterated expectations, can be expressed as $J^{*,d} = \min_{\vu = \gamma^d(\vy)}\mathbb{E}\left[\|\vm-\vu\|^2 \middle\vert \vy\right]$. Hence, a minimizer policy of the decoder is $\vu = \gamma^{*,d}(\vy) = \mathbb{E} \left[\vm \middle\vert \vy\right]$. Since both $\vm$ and $\vy$ are Gaussian, then the optimal decoder is
		\begin{align}
		E &\left[\vm \middle\vert \vy\right] = \mathbb{E}[\vm] + \Sigma_{\vM\vY}\Sigma_{\vY\vY}^{-1}(\vy-\mathbb{E}[\vy]) \nn\\
		&= \Sigma_{\vM}A^T \left(A\Sigma_{\vM}A^T + \Sigma_{\vW}\right)^{-1}(\vy-\vC)\,.
		\label{eq:decoderMultiPolicy}
		\end{align}
		\item Let the affine decoding policy be $\vu =\gamma^d(\vy) = K\vy+\vL$ where $K$ is an $n \times n$ matrix and $\vL$ is an $n \times 1$ vector. Then $\vu = K\vy+\vL = K(\vx+\vw)+\vL=K\gamma^e(\vm)+K\vw+\vL$. By using the completion of the squares method, the optimal cost is
		\begin{align*}
		J^{*,e} &= \min_{\vx = \gamma^e(\vm)} \mathbb{E} \left[\|\vm-\vu-\vb\|^2+\lambda \|x\|^2\right] \nn\\
		&= \mathbb{E} \left[\min_{\vx = \gamma^e(\vm)} \mathbb{E} \left[\|\vm-\vu-\vb\|^2+\lambda \|x\|^2 \Big|\vm \right] \right] \nn\\
		& = \mathbb{E} \Bigg[ \min_{\gamma^e(\vm)} \mathbb{E} \Big[\Big((K^TK+\lambda I)\gamma^e(\vm) - K^T\nn\\
		&\qquad\times(\vm-\vL-\vb)\Big)^T \Big(K^TK+\lambda I\Big)^{-1} \nn\\
		& \qquad\times\Big((K^TK+\lambda I)\gamma^e(\vm)- K^T(\vm-\vL-\vb)\Big) \nn\\
		&\quad + \Big(\vm-\vL-\vb\Big)^T \Big(I-K(K^TK+\lambda I)^{-1}K^T\Big)\nn\\
		& \qquad \times\Big(\vm-\vL-\vb\Big)\Big|\vm\Big]\Bigg]  + \mathbb{E} \Big[\vw^TK^TK\vw\Big]\,.
		\end{align*}
		Hence, the optimal $\gamma^e(m)$ can be chosen as follows:
		\begin{align}
		\gamma^{*,e} (\vm) &= A\vm + \vC = \Big(K^TK+\lambda I \Big)^{-1}K^T\nn\\
		& \qquad\qquad\qquad\qquad \times\Big(\vm-\vL-\vb\Big)\,.
		\label{eq:encoderMultiPolicy}
		\end{align}
		\item We have $K=\Sigma_{\vM}A^T \left(A\Sigma_{\vM}A^T + \Sigma_{\vW}\right)^{-1}$ and $A = \Big(K^TK+\lambda I \Big)^{-1}K^T$ from \eqref{eq:decoderMultiPolicy} and \eqref{eq:encoderMultiPolicy}. By combining these, $A = T(A)= \left(FF^T+\lambda I\right)^{-1}F$ can be obtained.
		\item Since $FF^T$ is a real and symmetric matrix, it is diagonalizable and can be written as $FF^T=Q \Upsilon Q^{-1}$ for a diagonal $\Upsilon$. Now consider $\|T(A)\|_F$ where $\|\cdot\|_F$ denotes the Frobenius norm:
		\begin{align}
		\|T(A)\|_F &= \tr \Bigg(\Big(\left(FF^T+\lambda I\right)^{-1}F\Big)\nn\\
		&\qquad\qquad\qquad \times\Big(\left(FF^T+\lambda I\right)^{-1}F\Big)^T\Bigg) \nn\\
		&= \tr \Bigg((\Upsilon+\lambda I)^{-1}\Upsilon (\Upsilon+\lambda I)^{-1}\Bigg) \nn\\
		&= \sum\limits_{i=1}^{n} {\upsilon_i \over (\upsilon_i+\lambda)^2}\,,
		\label{eq:multiDMappingBound}
		\end{align}
		where $\upsilon_i, i=1,\ldots,n$ are the eigenvalues of $FF^T$ and since $FF^T$ is positive semi-definite, all these eigenvalues are nonnegative. Since $\lambda>0$, we observe the following:
		\begin{align*}
		\upsilon_i \in [0,1] &\Rightarrow {\upsilon_i \over (\upsilon_i+\lambda)^2} < {1 \over \lambda^2} \,,\nn\\
		\upsilon_i \in (1,\infty) &\Rightarrow {\upsilon_i \over (\upsilon_i+\lambda)^2} < {\upsilon_i \over \upsilon_i^2} = {1 \over \upsilon_i} < 1\,.
		\end{align*}
		Hence, ${\upsilon_i / (\upsilon_i+\lambda)^2} < \max(1,{1 / \lambda^2})$ always holds. Then, by \eqref{eq:multiDMappingBound}, we have $\|T(A)\|_F < n\max(1,{1 / \lambda^2})$, which implies that $T(A)$ can be viewed as a continuous function mapping the compact convex set $\|A\|_F \in [0, n\max(1,{1 / \lambda^2})]$ to itself. Therefore, by Brouwer's fixed point theorem \cite{InfiniteDimensionalAnalysis}, there exists $A=T(A)$. \hfill$\IEEEQEDhere$
	\end{enumerate}
\end{IEEEproof}
We note, however, that there always exist a non-informative equilibrium (see Proposition~\ref{prop:nonInEq}, which also applies to the signaling game discussed in this section). However, there exist games with informative affine equilibria as we state in the following (see Theorem~\ref{thm:multiInformativeEq}).
\begin{proposition}
If either $\lambda$ or $\Sigma_{\vW}$ is zero, an informative affine equilibrium exists only if $\lambda$, $\Sigma_{\vW}$ and $\vb$ are all zero.
\end{proposition}
\begin{IEEEproof}
Note that, from \eqref{eq:decoderMultiPolicy} and \eqref{eq:encoderMultiPolicy}, we have $A=\Big(K^TK+\lambda I \Big)^{-1}K^T$, $K=\Sigma_{\vM}A^T \left(A\Sigma_{\vM}A^T + \Sigma_{\vW}\right)^{-1}$, $\vL=-K\vC$ and $A\vb=(AK-I)\vC$. From these equalities, we can analyze the equilibrium as in the scalar case:
\begin{enumerate}
	\item when $\lambda = 0$ and the noise is zero ($\Sigma_{\vW}=0$), then $A=K^{-1}$ and $K=A^{-1}$ are obtained. Then $A\vb=(AK-I)\vC=0$, thus the consistency of the equalities can be satisfied if only if $\vb=0$. Hence, if $\vb \neq 0$, there cannot exist an informative affine equilibrium. Recall that in the multi-dimensional noiseless cheap talk, the linearity of the equilibrium is shown for the uniform source; here the source is Gaussian.
	\item when $\lambda=0$, then $A=K^{-1}$ and $A\Sigma_{\vM}A^T + \Sigma_{\vW}=K^{-1}\Sigma_{\vM}A^T$ are obtained. There does not exist a solution to \eqref{eq:fixedPointEq} unless the noise is zero ($\Sigma_{\vW}=0$). Even when \eqref{eq:fixedPointEq} has a fixed point solution $A$, \eqref{eq:linearPolicyEq} and \eqref{eq:decoderPolicyEq} cannot hold together unless $\vb=0$.
	\item when the noise is zero ($\Sigma_{\vW}=0$), then $K=A^{-1}$ and $K^TK+\lambda I = K^TA^{-1}$ are obtained. There does not exist a solution to \eqref{eq:fixedPointEq} unless $\lambda=0$. Even when \eqref{eq:fixedPointEq} has a fixed point solution $A$,  \eqref{eq:linearPolicyEq} and \eqref{eq:decoderPolicyEq} cannot hold together unless $\vb=0$.\hfill$\IEEEQEDhere$
\end{enumerate} 
\end{IEEEproof}
\begin{remark}
In the multi-dimensional case, fixed points may not be unique: with $\lambda=1.0311$ and
\begin{align*}
\Sigma_{\vM}= \begin{bmatrix}
1.6421  &  0.1299  &  0.5713  &  0.2305 \\
0.1299  &  1.4803  &  0.6810  &  0.4749 \\
0.5713  &  0.6810  &  1.7312  &  0.4292 \\
0.2305  &  0.4749  &  0.4292  &  1.3515
\end{bmatrix} \,,\\
\Sigma_{\vW}= \begin{bmatrix}
1.2742  &  0.1868  &  0.2318  &  0.0559  \\
0.1868  &  1.8266  &  0.5955  &  0.3091  \\
0.2318  &  0.5955  &  1.2377  &  0.4951  \\
0.0559  &  0.3091  &  0.4951  &  1.5336
\end{bmatrix} \,,
\end{align*}
we can obtain two fixed points with different absolute-valued elements as follows (recall that if $A$ is a fixed point, $-A$ is also a fixed point):
\begin{align*}
A = \begin{bmatrix}
-0.1543 & 0.1762  &	0.0606	& 0.1117 \\	
0.1602	& 0.0159  & 0.1036	& 0.0279 \\	
-0.2000	& -0.1879 &	-0.2700	& -0.1565 \\	
0.0603	& 0.1052  &	0.1221	& 0.0824 	
\end{bmatrix} \,,\\
A = \begin{bmatrix}
-0.2431 & 0.0738 & -0.0752 & 0.0285	 \\
0.0293 & -0.1351 & -0.0966 & -0.0948 \\	
0.1520 & 0.2181 & 0.2682 & 0.1735	 \\
-0.1003 & -0.0801 & -0.1236 & -0.0683	
\end{bmatrix} \,.
\end{align*}
\end{remark}
\begin{thm}\label{thm:multiInformativeEq}
	Let source $\vM$ be a zero-mean $n$-dimensional Gaussian random variable with covariance matrix $\Sigma_{\vM}=diag\{\sigma_{m_1}^2,\ldots,\sigma_{m_n}^2\}$ where $diag$ indicates a diagonal matrix, and noise $\vW$ be a zero-mean $n$-dimensional Gaussian random variable with covariance matrix $\Sigma_{\vW}=diag\{\sigma_{w_1}^2,\ldots,\sigma_{w_n}^2\}$. Then an informative affine equilibrium exists if $\lambda<\max\{{\sigma_{m_1}^2\over\sigma_{w_1}^2},\ldots,{\sigma_{m_n}^2\over\sigma_{w_n}^2}\}$.
	\label{thm:multiIndependent}
\end{thm}
\begin{IEEEproof}
		Since the source components are independent and the noise components are independent, the $n$-dimensional noisy signaling game problem turns into $n$ independent scalar noisy signaling game problems as follows: 
		\begin{enumerate}
			\item If the decoder uses the channels independently; i.e., $u_i=\gamma_i^d(y_i)$ for $i=1,\ldots,n$, then the optimal cost of the encoder will be
			\begin{align*}
			& J^{*,e} = \min_{\vx = \gamma^e(\vm)} \mathbb{E} \left[\|\vm-\vu-\vb\|^2+\lambda \|x\|^2\right] \nn\\
			&\quad= \min_{\vx = \gamma^e(\vm)} \sum\limits_{i=1}^{n} \mathbb{E} [(m_i-\gamma_i^d(y_i)-b_i)^2+\lambda x_i^2] \nn\\
			&\quad= \sum\limits_{i=1}^{n} \min_{\vx = \gamma^e(\vm)} \mathbb{E} [(m_i-\gamma_i^d(y_i)-b_i)^2+\lambda x_i^2]\,.
			\end{align*}
			Since, $y_i=x_i+w_i$ for each $i=1,\ldots,n$, the optimal encoder also uses the channels independently; i.e., $x_i=\gamma_i^e(m_i)$ for $i=1,\ldots,n$.
			\item Similarly, if the encoder uses the channels independently; i.e., $x_i=\gamma_i^e(m_i)$ for $i=1,\ldots,n$, then the optimal cost of the decoder will be
			\begin{align*}
			J^{*,d} &= \min_{\vu = \gamma^d(\vy)} \mathbb{E} \left[\|\vm-\vu\|^2\right]\nn\\
			&= \sum\limits_{i=1}^{n} \min_{\vu = \gamma^d(\vy)} \mathbb{E} [(m_i-u_i)^2] \,.
			\end{align*}
			Since, $y_i=\gamma_i^e(m_i)+w_i$ for each $i=1,\ldots,n$, the optimal decoder will also use channels independently; i.e., $u_i=\gamma_i^d(y_i)$ for $i=1,\ldots,n$.
		\end{enumerate}
 Thus, we have, in each dimension $i$ ($i=1,\ldots,n$); 
	\begin{itemize}
		\item the source $M_i$ is a zero-mean Gaussian with variance $\sigma_{m_i}^2$,
		\item the channel has the Gaussian noise $W_i$ with zero-mean and variance $\sigma_{w_i}^2$,
		\item the encoder's goal is to find the optimal policy which minimizes its cost $\min_{x_i = \gamma^e(m_i)} \mathbb{E} [(m_i-u_i-b_i)^2+\lambda x_i^2]$,
		\item the decoder's goal is to find the optimal policy which minimizes its cost $\min_{u_i = \gamma^d(y_i)} \mathbb{E} [(m_i-u_i)^2]$.
	\end{itemize}
	For each dimension, the informative affine equilibrium exists if $\lambda<\sigma_{m_i}^2/\sigma_{w_i}^2$. For the multidimensional setup, the existence of the informative equilibrium in at least one dimension implies the existence of the informative equilibrium for the whole sytem. Hence, it is sufficient that the inequality $\lambda<\sigma_{m_i}^2/\sigma_{w_i}^2$ is valid for at least one dimension. As a result, the condition for the existence of the informative affine equilibrium becomes $\lambda<\max\{{\sigma_{m_1}^2\over\sigma_{w_1}^2},\ldots,{\sigma_{m_n}^2\over\sigma_{w_n}^2}\}$.
\end{IEEEproof}

Note that, from \eqref{eq:decoderMultiPolicy} and \eqref{eq:encoderMultiPolicy}, by assuming $|A| \neq 0$, we have $\lambda A\Sigma_{\vM}A^T = K^TK\Sigma_{\vW}$ 
which is equivalent to 
\begin{align}
\lambda (A^T)^{-1}\Sigma_{\vM}A^T = (K^TK+\lambda I)(K^TK+\lambda I)\Sigma_{\vW}\,.
\label{eq:multiShortEq}
\end{align}  
\begin{remark}
Assuming all channels are informative, i.e., $|A| \neq 0$, we make the following observations.
\begin{itemize}
	\item If the source is i.i.d.; i.e., $\Sigma_{\vM}=\sigma_m^2 I$, then \eqref{eq:multiShortEq} becomes
	\begin{align*}
	& \lambda (A^T)^{-1}\sigma_m^2 I A^T = (K^TK+\lambda I)(K^TK+\lambda I)\Sigma_{\vW} \nn\\
	&\Rightarrow \lambda\sigma_m^2(\Sigma_{\vW})^{-1} = (K^TK+\lambda I)(K^TK+\lambda I) \nn\\
	&\Rightarrow \lambda\sigma_m^2(\Sigma_{\vW})^{-1} \geq \lambda^2 I\nn\\
	& \Rightarrow \lambda I \leq \sigma_m^2(\Sigma_{\vW})^{-1}\,.
	\end{align*} 
		This result implies that $\lambda$ must satisfy the inequality $\lambda I \leq \sigma_m^2(\Sigma_{\vW})^{-1}$ for the i.i.d. source; otherwise, there must be at least one non-informative channel; i.e., $|A|$ must be $0$.
	\item If the channel noise is i.i.d.; i.e., $\Sigma_{\vW}=\sigma_w^2 I$, (since $\Sigma_{\vM}$ is real-symmetric, it has the eigenvalue decomposition as $\Sigma_{\vM}=Q\Lambda Q^T$), then \eqref{eq:multiShortEq} becomes
	\begin{align*}
	& \lambda (A^T)^{-1}\Sigma_{\vM}A^T = (K^TK+\lambda I)(K^TK+\lambda I)\sigma_w^2 I \nn\\
	&\Rightarrow {\lambda\over\sigma_w^2} (A^T)^{-1}Q\Lambda Q^T A^T = (K^TK+\lambda I)(K^TK+\lambda I)\nn\\
	&\Rightarrow (A^T)^{-1}Q\Lambda Q^T A^T \geq \lambda\sigma_w^2\,.
	\end{align*} 
	This result implies that for each eigenvalue $\lambda_{\vM}$ of $\Sigma_{\vM}$, $\lambda$ must satisfy $\lambda\leq\lambda_{\vM} /\sigma_w^2$ for the i.i.d. channel noise; otherwise, there must be at least one non-informative channel; i.e., $|A|$ must be $0$.
	\item For the general case, recall the Minkowski determinant theorem, $|A+B|^{1/n}\geq|A|^{1/n}+|B|^{1/n}$ which holds for any non-negative $n\times n$ Hermitian matrix $A$ and $B$. This implies $|A+B|\geq|A|+|B|$. By using this inequality and \eqref{eq:encoderMultiPolicy},
	\[ |A| = {|K|\over|K^TK+\lambda I|}\leq{|K|\over|K|^2+\lambda^n}\,.\]
	Assuming $|A|\neq0$, recall the equality $\lambda A\Sigma_{\vM}A^T = K^TK\Sigma_{\vW}$. Taking the determinant of both sides, 
	\begin{align*}
	|K|^2 |\Sigma_{\vW}| &= \lambda^n |A|^2 |\Sigma_{\vM}| \leq \lambda^n \left({|K|\over|K|^2+\lambda^n}\right)^2 |\Sigma_{\vM}| \nn\\
	&\leq \lambda^n{|K|^2\over\lambda^{2n}}|\Sigma_{\vM}| \Rightarrow\lambda^n\leq{|\Sigma_{\vM}|\over|\Sigma_{\vW}|}\,.
	\end{align*}
	The result can be interpreted as follows: If $\lambda>\left({|\Sigma_{\vM}|\over|\Sigma_{\vW}|}\right)^{1/n}$, then $|A|=|K|=0$ in the equilibrium; i.e., there must be at least one non-informative channel.
\end{itemize}
\end{remark}
 
\subsection{An information theoretic lower bound on the encoder cost and the existence of informative equilibria }
We end the section with an information theoretic lower bound for the encoder cost; this serves us also to obtain condition for the existence of an informative equilibrium.
Let $\ve=\vm-\vu=\vm-\mathbb{E}[\vm|\vy]$, then we have $\Sigma_{\ve}=\mathbb{E}[\ve\ve^T]=\mathbb{E}[(\vm-\mathbb{E}[\vm|\vy])(\vm-\mathbb{E}[\vm|\vy])^T]$. From the information theoretic inequalities;
\begin{align*}
I(\vm;\vy) &= h(\vm)-h(\vm|\vy) = h(\vm)-h(\vm-\mathbb{E}[\vm|\vy]|\vy) \nn\\
&\geq h(\vm)-h(\vm-\mathbb{E}[\vm|\vy]) \nn\\
&\geq h(\vm)-{1\over 2}\log_2((2\pi e)^n |\Sigma_{\ve}|) \nn\\
& = {1\over 2}\log_2((2\pi e)^n |\Sigma_{\vm}|)-{1\over 2}\log_2((2\pi e)^n |\Sigma_{\ve}|) \nn\\
& = {1\over 2}\log_2(|\Sigma_{\vm}|/|\Sigma_{\ve}|)\,.
\end{align*} 
Also from the rate-distortion theorem, the data processing theorem and the channel capacity theorem:
\begin{align*}
R(D)&\leq \min_{f(\vu|\vm):\mathbb{E}[\|\vm-\vu\|^2]\leq D}I(\vm;\vu) \leq I(\vm;\vu) \nn\\
& \leq I(\vx;\vy) \leq \max_{f(\vx):\mathbb{E}[\|\vx\|^2]\leq P}I(\vx;\vy) \leq C(P)\,.
\end{align*}
If we combine these, we obtain the following:
\begin{align}
|\Sigma_{\ve}| &\geq |\Sigma_{\vm}|2^{-2R(D)}\geq |\Sigma_{\vm}|2^{-2I(\vm;\vu)} \nn\\
&\geq |\Sigma_{\vm}|2^{-2I(\vx;\vy)}\geq |\Sigma_{\vm}|2^{-2C(P)}\,.
\label{eq:rateToCapacity}
\end{align}
Now consider the following:
\begin{align}
& \mathbb{E}[\|\vm-\vu\|^2]=\mathbb{E}[\|\ve\|^2]=\tr\Sigma_{\ve}\overset{(a)}{\geq} n\Big(\prod_{i=1}^{n}\Sigma_{\ve}(i,i)\Big)^{1/n}\nn\\
&\qquad\qquad\overset{(b)}{\geq}n\Big(|\Sigma_{\ve}|\Big)^{1/n}\overset{(c)}{\geq}n\Big(|\Sigma_{\vm}|2^{-2C(P)}\Big)^{1/n}\,.
\label{eq:boundWithCapacity}
\end{align}
Here, (a) follows from the inequality for the arithmetic and geometric mean where $\Sigma_{\ve}(i,i)$ stands for $i$th diagonal element of $\Sigma_{\ve}$, (b) follows from the Hadamard inequality (since $\Sigma_{\ve}$ is a positive semi-definite matrix), and (c) follows from \eqref{eq:rateToCapacity}.
Now we will rewrite \cite[Eq. (9.166)]{thomasCover} which presents the capacity of the additive colored Gaussian noise channel with typo corrected:
\begin{align*}
C(P)={1\over n}\sum\limits_{i=1}^{n}{1\over 2}\log_2\Big(1+{\max(\nu-\lambda_i,0)\over\lambda_i}\Big)\,,
\end{align*}
where $P=\mathbb{E}[\|\vx\|^2]$, $\lambda_1, \lambda_2,\ldots,\lambda_n$ are the eigenvalues of $\Sigma_{\vw}$ and $\nu$ is chosen so that $\sum_{i=1}^{n}\max(\nu-\lambda_i,0)=nP$. Then we can obtain the following:
\begin{align}
2^{-2C(P)}&=2^{-2{1\over n}\sum\limits_{i=1}^{n}{1\over 2}\log_2\Big(1+{\max(\nu-\lambda_i,0)\over\lambda_i}\Big)}\nn\\
&=\prod_{i=1}^{n}\Big(1+{\max(\nu-\lambda_i,0)\over\lambda_i}\Big)^{-1/n}\nn\\
&=\prod_{i=1}^{n}\Big({\max(\nu,\lambda_i)\over\lambda_i}\Big)^{-1/n}\nn\\
&={(\prod_{i=1}^{n}\lambda_i)^{1/n}\over(\prod_{i=1}^{n} \max(\nu,\lambda_i))^{1/n}}\nn\\
&\overset{(a)}{\geq}{(|\Sigma_{\vw}|)^{1/n}\over(P+\sum_{i=1}^{n} {\lambda_i\over n})}\nn\\
&=\Big(|\Sigma_{\vw}|\Big)^{1/n}\Big(P+{\tr\Sigma_{\vw}\over n}\Big)^{-1}\,.
\label{eq:capacityBound}
\end{align}
Here, (a) holds, since our assumption $\sum_{i=1}^{n}\max(\nu-\lambda_i,0)=nP$ implies $\sum_{i=1}^{n}\max(\nu,\lambda_i)=nP+\sum_{i=1}^{n}\lambda_i$ and $(\prod_{i=1}^{n} \max(\nu,\lambda_i))^{1/n}\leq\sum_{i=1}^{n} \max(\nu,\lambda_i)/n=P+\sum_{i=1}^{n} \lambda_i/ n$ holds by the inequality for the arithmetic and geometric mean. If we insert \eqref{eq:capacityBound} to \eqref{eq:boundWithCapacity},
\begin{align}
& \mathbb{E}[\|\vm-\vu\|^2] \geq n\Bigg(|\Sigma_{\vm}|\Big(|\Sigma_{\vw}|\Big)^{1/n}\Big(P+{\tr\Sigma_{\vw}\over n}\Big)^{-1}\Bigg)^{1/n} \nn\\
&\qquad\quad= n(|\Sigma_{\vm}|)^{1/n}(|\Sigma_{\vw}|)^{1/n^2}\Big(P+{\tr\Sigma_{\vw}\over n}\Big)^{-1/n}\,.
\label{eq:errorBound}
\end{align}
The encoder costs reduces to $J^e = \mathbb{E}[\|\vm-\vu\|^2]+\lambda \mathbb{E}[\|\vx\|^2]+\|\vb\|^2$ since the decoder always chooses $\vu=\mathbb{E}[\vm|\vy]$. Then, by \eqref{eq:errorBound},
\begin{align}
J^e &= \|\vb\|^2 + \lambda \mathbb{E}[\|\vx\|^2] + \mathbb{E}[\|\vm-\vu\|^2] \nn\\
& \geq \|\vb\|^2 + \lambda P + n(|\Sigma_{\vm}|)^{1/n}(|\Sigma_{\vw}|)^{1/n^2}\Big(P+{\tr\Sigma_{\vw}\over n}\Big)^{-1/n}\,.
\label{eq:itAsBound}
\end{align}
The minimizer of this function can be found by the local perturbation condition:
\begin{align*}
\lambda-&(|\Sigma_{\vm}|)^{1/n}(|\Sigma_{\vw}|)^{1/n^2}\Big(P+{\tr\Sigma_{\vw}\over n}\Big)^{-{1\over n}-1}=0 \nn\\
\Rightarrow\lambda&=(|\Sigma_{\vm}|)^{1/n}(|\Sigma_{\vw}|)^{1/n^2}\Big(P+{\tr\Sigma_{\vw}\over n}\Big)^{-{1\over n}-1} \nn\\
&\overset{(a)}{\leq}(|\Sigma_{\vm}|)^{1/n}(|\Sigma_{\vw}|)^{1/n^2}\Big((|\Sigma_{\vw}|)^{1/n}\Big)^{-{1\over n}-1}\nn\\
&=(|\Sigma_{\vm}|)^{1/n}(|\Sigma_{\vw}|)^{-1/n}\,.
\end{align*}
Here, (a) follows from the nonnegativeness of $P$ and the inequality for the arithmetic and geometric mean and the Hadamard inequality, similar to \eqref{eq:boundWithCapacity}. Hence, if $\lambda<(|\Sigma_{\vm}|)^{1/n}(|\Sigma_{\vw}|)^{-1/n}$, the lower bound is minimized at a nonzero $P$ value, but if  $\lambda\geq(|\Sigma_{\vm}|)^{1/n}(|\Sigma_{\vw}|)^{-1/n}$, the minimizer $P$ becomes zero. 
Finally, if channels and source are assumed to be i.i.d.; i.e., $\Sigma_{\vm}=\sigma_m^2 \mathbb{I}$ and $\Sigma_{\vw}=\sigma_w^2 \mathbb{I}$ where $\mathbb{I}$ is $n\times n$ identity matrix, and the encoder and the decoder use linear policies, then \eqref{eq:itAsBound} becomes tight and can be interpreted as follows: If $\lambda>(|\Sigma_{\vm}|)^{1/n}(|\Sigma_{\vw}|)^{-1/n}=\sigma_m^2/\sigma_w^2$, then \eqref{eq:itAsBound} is minimized at $P=0$; that is, the encoder does not signal any output. Hence, the encoder engages in an non-informative equilibrium and the minimum cost becomes $E\left[\|\vm\|^2\right]+\|\vb\|^2$ at this non-informative equilibrium. Recall that this is analogous to the analysis in the scalar setup \eqref{eq:itInequality}.

\section{Concluding Remarks}

For a strategic information transmission problem under quadratic criteria with a non-zero bias term leading to a mismatch in the encoder and the decoder objective functions, Nash and Stackelberg equilibria have been investigated in a number of setups. It has been proven that for any scalar source, the quantized nature of Nash equilibrium policies hold, whereas all Stackelberg equilibrium policies are fully informative. Further, it has been shown that the Nash equilibrium policies may be non-discrete and even linear for a multi-dimensional cheap talk problem, unlike the scalar case. The additive noisy channel setup with Gaussian statistics has also been studied, such a case leads to a signaling game due to the communication constraints in the transmission. Conditions for the existence of affine Nash equilibrium policies as well as general informative Nash equilibria are presented for both the scalar and multi-dimensional setups. Lastly, we proved that the only equilibrium in the Stackelberg noisy setup is the linear equilibrium. Our findings provide further conditions on when affine policies may be optimal in decentralized multi-criteria control problems and lead to conditions for the presence of active information transmission in strategic environments. Recently, we have extended our analysis in this paper to study dynamic signaling games \cite{ISIT2016dynamicSignalingGame}.

\bibliographystyle{IEEEtran}
\bibliography{SerdarBibliography}

\end{document}